# The Laplace representation of $1/(\sin(\pi s/4)\cdot 2\xi(\frac{1}{2} + s))$.

## On the Riemann zeta-function, Part II.

By Anthony Csizmazia

E-mail: apcsi2000@yahoo.com

**Abstract**

An odd meromorphic function f(s) is constructed from the Riemann zeta-function evaluated at one-half plus s. We determine the two-sided Laplace transform representation of f(s) on open vertical strips, V'(4w), disjoint from the (translated) critical strip. V'(4w) consists of all s with real part, Re(s), of absolute value greater than one-half and Re(s) between successive poles 4w, 4(w + 1) of f(s), with w an integer. The corresponding Laplace density is related to confluent hypergeometric functions. That density is shown to be positive for nonzero w other than -1. Those results are obtained without relying on any unproven hypothesis. They are used together with the Riemann hypothesis and hypotheses advanced by the author to obtain conditional results concerning the zeta-function. Those results are presented in Part I. Their proofs are derived in Parts III-V. A metric geometry expression of the positivity of the Laplace densities arising is established in Part VI.

**Keywords:** Riemann zeta-function; Laplace transform; analytic / entire / meromorphic / function; confluent hypergeometric function; positive definite function; analytic characteristic function.

**MSC (Mathematics Subject Classification).** 11Mxx Zeta and L-functions: analytic theory. 11M06 $\zeta(s)$ and $L(s, \chi)$. 30xx Functions of a complex variable. 44A10 Laplace transform. 42xx Fourier analysis. 42A82 Positive definite functions. 60E10 Characteristic functions; other transforms. 33C15 Confluent hypergeometric functions, Whittaker functions, ${}_1F_1$.

**Table of contents.**

Properties of $\zeta(s)$. Property of $\zeta(s)$. Integrability of $|m(x + it, \beta)|$ in t. The Laplace representation of $m(z, \beta)$. The Laplace representation of $\pi/\sin(\pi s)$. Translation principal.

**§2 The Laplace representation on $V_0'$ of $f(s) := 1/(b(s)\zeta(\frac{1}{2} + s))$ via that of $1/b(s)$.**

**§3 Strategy for determining the Laplace representation on $V_0'$ of $1/b(s)$.**

Divergence of the formal partial fraction expansion of $1/b(s)$. Integrability of $|f(s)|$ on vertical lines outside of the critical strip. Bound for $1 / |\zeta(z)|$. Integrability of $|N(z, \frac{1}{4})|$ on vertical lines. Determination of the Laplace representations of $f_0(s)$ and $1/b(s)$ on $V_0'$ from that of $F(z, \frac{1}{4})$ on $V(\frac{1}{4}, 2)$. Translation relation for $F(z, \beta)$. Translation relation for $f(z)$. The Mellin transform representation of $j(u, m)$. The Mellin transform representation of $1/(z)_m$. The Mellin transform representations of $1/(S + 2m - \frac{1}{2})$, $j(\frac{1}{2} + S, m)$ and $f_0(S)$. The Mellin transform representation of $f(s)$, for s on $V_0' + 4w$. Positivity. Determination of the Laplace representation of $f(s)$ on $V_{4w}'$. Determination of the Laplace representation of $F(z, \beta)$ relative to z on $V_{4w}$. Integrability of $|F(z, \beta)|$ on vertical lines. Divergence of the formal partial fraction expansion of $F(z, \beta)$. Splitting $F(u, \beta, 1)$. Trigonometric identity. The splitting of $F(u, \beta, 1)$ via $E(u, \beta, 1)$ and $E(u, \beta, 2)$. Trigonometric asymptotics. The asymptotic behavior of $|E(u, \beta, 1)|$ and of $|E(u, \beta, 2)|$. Integrability. The Mellin transform representation of $E(u, \beta, 2)$. The Mellin transform representation of $E(u, \beta, 1)$. Definition and properties of $W(z, \beta)$. Translation relation for $W(z, \beta)$. Mellin transform representation of $((\pi/2)/\cos((\pi/2)(u + \beta)))\cdot\Gamma(u)$. Definition and properties of $B_0(z, \beta)$. Mellin transform representation of $E(u, \beta, 1)$. The Mellin transform representation of $F(u, \beta, 1)$. Definition and properties of $M(z, \beta)$. Mellin transform representation of $F(u, \beta, 1)$. Definition and properties of $I(p, z, u)$. The determination of $W(z, 1 + \beta)$ from $I(p, z, u)$. The determination of $I(p, z, u)$ from $I(p, z/u)$. The determination of $W(z, 1 + \beta)$ from $I(-\beta, \pm iz)$.

**§4 The incomplete gamma functions and confluent hypergeometric functions**

Relation between $I(p, z)$ and $\Gamma(1 - p, z)$. The confluent hypergeometric functions $M(a, B, z)$ and $U(a, a, z)$. Kummer's equation. Relation between $\varphi(1 + \beta, z)$ and $\gamma(\beta, z, *)$. The Laplace representation of $\varphi(1 + \beta, z)$. Relation between $\varphi(1 + \beta, z)$ and $\varphi(\beta, z)$. Relation between $I(-\beta, z)$ and $\varphi(1 + \beta, z)$.

**§5 $H(z, \beta)$ and the Mellin transform representation of $F(z, \beta)$.**

The determination of $W(z, 1 + \beta)$ from $H(z, \beta)$. The determination of $M(z, \beta)$

from W(z, 1 + β). The determination of M(z, β) from H(z, β). Integrability of |F(z, β)| on vertical lines. Representation of H(z, β). Uniform boundedness. Positivity. Translation relation for H(z, β). The determination of H(z, β) from W(z, β - 1). The representation of W(z, β) via R(z, β). The boundedness of H(z, β). Mellin transform representation of F(z, β). Positivity.

**§6 The Mellin transform representation of f(s, β) := 1/(sin(πs/4)·2ξ(2β + s)).**

Relations used to determine the Mellin transform representation of f(s, β). Integrability of 1/|b(s, β)| on vertical lines. The Mellin transform representation of $f_0$(s, β). The Mellin transform representation of 1/b(s, β). The Mellin transform representation of f(s, β). The determination of $P_{4w}$ from $T_0$. The order of $P_0$(r, β) in r. Results when β = ¼ . Main unconditional theorem (1), (4) (i). Positivity when β and w are nonnegative. Results when β = ¼ . Main unconditional theorem (4) (i), (ii). Metric norms and analytic characteristic functions. Metric result when β = ¼ .

**References**

**Index of abbreviations.**

**§1**
Corollary 1.1 (Analytic characteristic function) ACF.

**Index of symbols**

(Complex plane) C. (Real line) **R**, $\mathbb{R}$.

**§1**
n(s, α), f(s, α). m(z, β). $Q_k$(z). M(x). q(u, β), $l_k$(y, β).

**§2**
P(z), g(r, j, q), g(j, q)(r), g(r, j, q, p). $q_0$(T), $j_0$(T). $\theta_{T,p}$(r), $\omega_{T,p}$(r).

**§3**
$n_0$(s), $f_0$(s), $h^{<\alpha>}$. N(z, β), F(z, β). j(u, m), E(v, m) N(u, β, 1), F(u, β, 1). E(u, β, 1), E(u, β, 2), J(z). W(z, β), $B_0$(z, β), M(z, β). I(p, z, u), I(p, z).

**§4**
φ(B, z), γ(β, z, ∗), (confluent hypergeometric functions) M(a, B, z). (Kummer's operator) K(a, B). U(a, a, z)

**§5**

$H(z, \beta)$. $\Omega(r, \delta)$. $A(z, \beta, n)$, $R(z, \beta)$. $G(z, \beta, m)$, $H(z, \beta, 2w)$.

**§6**
$n_0(s, \beta)$, $f_0(s, \beta)$, $b(s, \beta)$, $T_0(z, \beta)$, $T_0(z, \beta, 4w)$. $n(s, \beta)$, $f(s, \beta)$, $c(z, \beta)$, $P_{4w}(z, \beta)$. $m(t)$, $d(t_1, t_2)$.

**The Laplace representation of $1/(\sin(\pi s/4)\cdot 2\xi(\frac{1}{2} + s))$.**

**On the Riemann zeta-function, Part II.**

**Review of elements of Part I.** See A. Csizmazia [5].

**§1 Definitions of $l(s)$, $a(s)$, $\xi(s)$, $n(s)$, $f(s)$, $b(s)$, $V(x_0, x_1)$, $V[x_0, x_1]$, $V(\varepsilon)$.**

Let s be complex. Define:

$$l(s) := \pi^{-s/2} s\Gamma(s/2) = \pi^{-s/2}\cdot 2\Gamma(1 + s/2),\ a(s) := l(s)(s-1),\ \xi(s) := (\tfrac{1}{2})a(s)\zeta(s),$$

with $\zeta(s)$ the Riemann zeta-function,

$$n(s) := \sin(\pi s/4)\cdot 2\xi(\tfrac{1}{2} + s),\ f(s) := 1/n(s) \text{ and } b(s) := \sin(\pi s/4)a(\tfrac{1}{2} + s).$$

Say $x_0 < x_1$. Let $V(x_0, x_1)$ be the open vertical strip of all s with $x_0 < \mathrm{Re}(s) < x_1$. Define $V[x_0, x_1]$ to be the closed strip of all s with $x_0 \le \mathrm{Re}(s) \le x_1$. Set $V(\varepsilon) := V(0, \varepsilon)$ for $\varepsilon > 0$.

§3

**Definitions of $V_u$, $V_u'$.** Say u is a multiple of four, $u = 4w$. Let $V_u := V(u, u+4)$. If $u \neq 0, -4$, set $V_u' = V_u$. Let $V_0' := V(\frac{1}{2}, 4)$. Take $V_{-4}' = -V_0'$.

**Definition of the Pochhammer symbol $(z)_n$.**
Say z is complex. Set $(z)_n := \Pi_{0\le k\le n-1}(z + k)$, if $n \ge 1$. Take $(z)_0 = 1$.

**Definitions of $\tilde{c}(4k)$, $c(z)$.** Let k be an integer $\ge 0$. Set

$$\tilde{c}(4k) := 1/(\pi^{3/4}\Gamma(5/4 + 2k)(2k - \tfrac{1}{4})\zeta(\tfrac{1}{2} + 4k)).$$

Define $c(z) := 1/n'(z)$, for z with $n'(z) \neq 0$. $c(4k) = \tilde{c}(4k)(-(\pi^2))^k$.

**Definition of $P_0(z)$.** Set $P_0(z) := (-1)\sum_{k\ge 1} \tilde{c}(4k)(-(z^2))^k$.

**Definition of the open disk $B(z, r)$.** Say $r \ge 0$. Take $B(z, r) := \{s: |s - z| \le r\}$.

**Definition of $P_{4w}(z)$.** Say w = 1,2,3,... Let

$$P_{4w}(z) := (-1)^w(P_0(z) + \sum_{1 \le k \le w} \tilde{\mathbf{c}}(4k)(-(z^2))^k) = (-1)^{w+1}\sum_{k \ge w+1} \tilde{\mathbf{c}}(4k)(-(z^2))^k.$$

**§4**

**(4.5)**

**Definition of an analytic characteristic function on a vertical strip.**
Say $-\infty \le w_0 < w_1 \le \infty$. j(s) is an "analytic characteristic function on $V(w_0, w_1)$" when j(s) is analytic on $V(w_0, w_1)$ and j(iz) is positive definite in z on the horizontal strip $-iV(w_0, w_1)$: $j(s) = \int_{\mathbf{R}} d(y)e^{sy}d(\mu)$, with μ a positive measure on the real line. Say $w_0 < x < w_1$ and t is real. Then j(x + it)/j(x) is the characteristic function relative to t of the probability measure $\mu_x(S) := \int_S d(y)(e^{xy}/j(s))d(\mu)$. See E. Lukacs [10].

**Definition of a meromorphic characteristic function on C.**
Let us say that h(z) is a "meromorphic characteristic function on C" when each of the following conditions hold.
(1) h(z) is the reciprocal of an entire function, j(z) say.
(2) All of the zeros of j(z) are on the union of the real axis and the imaginary axis.
(3) j(0) = 0. j(it) = 0 for at least one nonzero real t.
(4) The real zeros of j(z) are unbounded from above and also from below. Each real zero is simple.
(5) Let the successive real zeros be $w_k$ with k any integer, $w_0 = 0$ and $w_k < w_{k+1}$. Either h(z) or –h(z) is an analytic characteristic function on the open strip bounded by the vertical lines through $w_k$, $w_{k+1}$.

**Unconditional results. Outline of strategy and proofs.**

**Review** Part I, §1, Definitions; and §4, (4.5), Definition of an analytic characteristic function on a vertical strip, Definition of a meromorphic characteristic function on C.

**§1 The role of sin(πs/4) in f(s) := 1/(sin(πs/4)2ξ(½ + s)).**

We now elucidate the contribution of the factor sin(πs/4) of n(s) := sin(πs/4)2ξ(½ + s) toward obtaining that f(s) := 1/n(s) should be a meromorphic characteristic function on the complex plane C.

Let x, t be real.

**Claim 1.1**. *Say $0 < \varepsilon < \pi/4$. Fix x.*
(1) *$h(t) := 1/(e^{\varepsilon|t|}\xi(x + i\, t))$ is unbounded on the vertical line through x. Thus h(t) is not the Fourier transform of an absolutely integrable function on the real line.*
(2) *$\int_R |h(t)|^p\, dt$ is infinite, for any $p \geq 0$.*

**Claim 1.2**. *Say $0 < \varepsilon < \pi / 4$.*
(1) *Fix x. There exists $T = T(\varepsilon , x) \geq 0$ such that*
*(*): $1 / |\xi(x + i\, t)| > e^{\varepsilon|t|}$ if $|t| > T$.*
(2) *Say $1 \leq x_0 < x_1$. Then there exists a $T = T(\varepsilon, x_0, x_1) \geq 0$ such that for all x between $x_0$ and $x_1$ the inequality (*) of (1) holds.*

**Proof of Claims 1.1-1.2.** $\xi(s) := \pi^{-s/2}\Gamma(1 + (s/2))(s - 1)\zeta(s)$. The Stirling approximation for $\Gamma(z)$ (See Part I, §4, (4.1)) yields the following. Say $0 < \varepsilon < \pi/2$ and $x_0 < x_1$. Then there exists a $T = T(\varepsilon, x_0 , x_1) \geq 0$ such that $|\Gamma(x + i\, t)| < e^{-\varepsilon|t|}$ if x is between $x_0$, $x_1$ and $|t| \geq T$.

We now appeal to the following. See T. M. Apostol [3].
**Properties of ζ(s).**
(1) If $1 < \sigma \leq \mathrm{Re}(s)$, then $|\zeta(s)| \leq \zeta(\sigma)$.
(2) There exists $M > 0$ such that for $\sigma \geq 1$ and $|t| \geq e$, one has $|\zeta(\sigma + i\, t)| < M \log(|t|)$.
(3) If $\frac{1}{2} < \sigma \leq 1$ then $\zeta(\sigma + it) = O(|t|^{(1-\sigma)/2})$ .
(4) $\zeta(\frac{1}{2} + it) = O(|t|^{\frac{1}{2}})$ .

Then Claim 1.2 follows. Hence Claim 1.1 holds.

**Review** Part I §4, (4.5), Definitions of an analytic characteristic function on a vertical strip and of a meromorphic characteristic function on C. §7 Main conditional theorem (3). Definitions of ridge, groove functions.

Next it is revealed how to eliminate the problems expressed in Claims 1.1-1.2 by replacing $\xi(\frac{1}{2} + s)$ with $J(s) := j(s)\cdot 2\xi(\frac{1}{2} + s)$. The multiplier j(s) is to be such that J(s) is a minimum perturbation of $\xi(\frac{1}{2} + s)$ which yields that 1/J(s) is a meromorphic characteristic function on the complex plane. Obtaining that property for 1/J(s) must make a maximum reliance on the behavior of the zeta–function as opposed to that of j(s). It is to be anticipated that j(s) is optimal in a sense developed below in discussing the choice $j(s) := \sin(\pi s/4)$. Also J(s) should essentially conserve the relevant properties of $\xi(\frac{1}{2} + s)$.

$\theta(s) = \xi(½ + s)$ has the properties: $\theta(s)$ is entire, $\theta(s^*) = (\theta(s))^*$ and $\theta(-s) = \theta(\sigma s)$ for all s with $\sigma = 1$. Require that j(s) has those properties with $\sigma = 1$ or with $\sigma = -1$ in the last one. Then J(s) shares those properties, $\sigma$ being the same as for j(s).

Require that all of the roots z (with $Re(z) \geq 0$) of the entire function j(s) be real. Then the nonreal roots (respectively: real roots) of J(s) are precisely the roots of $\xi(½ + s)$ (respectively: j(s)) with the order of each root conserved.

Say $\varepsilon > 0$ and h(s) is an entire function on an open set containing the s with $0 \leq Re(s) < \varepsilon$. Assume 1/h(s) is an analytic characteristic function on $V(\varepsilon)$. Then h(s) is an extended groove function on $V(\varepsilon)$. Say t is real and $h(i{\cdot}t) = 0$. Then $h(0) = 0$. If $h'(0)$ is nonzero, then each root of h(s) on the imaginary axis is simple: $h'(i{\cdot}t) \neq 0$. Therefore require that 0 be a simple root of j(s). Thus j and J are odd.

Assume that each positive root x of j(s) is at least ½. In that case, if J(s) is an extended groove function on V(½), then RH and SZC hold.

Assume that all of the roots of j(s) are simple. That and SZC together imply that all roots of J(s) are simple. Then the formal partial fraction representation of 1 / J(s) has its simplest form. (See §5 of Part I.) That representation is used to derive the Laplace representation of 1/J(s) on V(½).

So consider the family of choices $j(s) := \sin(\alpha s)$ with $\alpha$ a nonzero real.

**Defininitions of n(s, α), f(s, α)**. Take $\alpha > 0$. Set

$$n(s, \alpha) := \sin(\alpha s)\cdot 2\xi(½ + s) \text{ and } f(s, \alpha) := 1/n(s, \alpha).$$

Assume $|t| \geq \delta > 0$. Then $1 / |\sin(x+it)| = 2\cdot e^{-|t|}(1+ e^{-2|t|}\varepsilon(x, t))$ with $-1 < \varepsilon(x, t) < 1 + 1/\delta$.

If $0 < \alpha < \pi/4$, then $n(s, \alpha)$ decays with exponential rapidity on any vertical line as $Im(s) \to \pm\infty$. Thus $f(s, \alpha)$ behaves like $1/\xi(s)$ does in Claim 1.1 above. So assume $\alpha \geq \pi/4$.

Say $x_0 < x < x_1$, t is real and $s = x + it$. Say |t| is large. $\sin(\alpha s)\cdot a(½ + s) \sim e^{(\alpha - (\pi/4))|t|} \cdot K_1(x)\cdot t^{7/4 + x/2} (1 + \varepsilon(x, t)/t)$, with $K_1(x)$ and $\varepsilon(x, t)$ as in §4.1 of Part I.

**Property of ζ(s).** If $Re(z) \geq \sigma > 1$, then $|\zeta(z)| \geq \zeta(2\sigma)/\zeta(\sigma)$. (This follows from the Euler factorization of $\zeta(z)$.) See T. M. Apostol [3].

Assume $\alpha > \pi/4$. Say $1 < x_0 \leq Re(s) \leq x_1$. Then $f(s, \alpha)$ converges to zero with uniform exponential rapidity as $Im(s) \to \pm\infty$.

Say $\alpha$ is nonzero. Then $\sin(\alpha s) = 0$ iff $s = k\pi/\alpha$ with k an integer.

Assume $\alpha > 0$. In that case: $n(x, \alpha)$ is nonzero for all x with $0 < x < ½$ iff $\alpha \leq 2\pi$.

Take $\alpha > 2\pi$. There is a unique integer $k_0$ with $k_0\pi/\alpha < ½ \leq (k_0 + 1)\pi/\alpha$ . Then $k_0 \geq 1$. Even if $n(x, \alpha)$ is an extended groove function on each $V(k\pi/\alpha, (k+1)\pi/\alpha)$ with $1 \leq k \leq k_0$ , it is still possible that RH fails for $\zeta(s)$ on the vertical lines through $½ + k\pi/\alpha$ for those k. That aspect worsens as $\alpha \to \infty$.

Assume $\alpha \geq \pi/4$. The choice $\alpha = \pi/4$ will prove to be optimal as expressed in the following lemma.

**Lemma 1.1** *Assume $(-1)^w f(s)$ is an analytic characteristic function on $V_{4w}$. Say $\alpha > \pi/4$, $4w \leq k\pi/\alpha$ and $(k+1)\pi/\alpha \leq 4(w+1)$. Then $(-1)^k f(s, \alpha)$ is an analytic characteristic function on $V(k\pi/\alpha, (k+1)\pi/\alpha)$.*

The previous lemma is a corollary of the results presented next.

**Defininition of m(z, β).** Let β, z be complex. When β is nonzero and z is not an integer, set $m(z, \beta) := (\pi/\beta)(\sin(\beta z)/\sin(\pi z))$. Then m has a continuous (analytic) extension also allowing $\beta = 0$ or $z = 0$. If z is not an integer, let $m(z, 0) := \pi z/\sin(\pi z)$. Set $m(0, \beta) :=1$ for any β. $m(\sigma_1 z, \sigma_2\beta) = m(z, \beta)$ for $\sigma_k = \pm 1$.

Assume β is nonzero and $-1 < \beta < 1$. Say $\delta > 0$. Let $s = x + it$ with x, t real and also $|t| \geq \delta$. Then $|\sin(\beta s)/\sin(s)| = e^{-(1 - |\beta|)|t|}(1+ e^{-2|\beta||t|}\varepsilon(x, t))$ with $|\varepsilon(x, t)|$ uniformly bounded, $|\varepsilon(x, t)| < K(\delta)$.

**Integrability of |m(x + it, β)| in t.** Say $-\pi < \beta < \pi$. Fix β. Then $|m(z, \beta)|$ is integrable in z on any vertical line not passing through a nonzero integer.

Fix $\alpha > \pi/4$. Say s is not a multiple of 4 or of $\pi/\alpha$. Then $f(s,\alpha) = (\pi/(4\alpha))m(z, \beta)\cdot f(s)$ with $z = (\alpha/\pi)s$, $\beta = \pi^2/(4\alpha)$ and $0 < \beta < \pi$. Assume $4w \leq k\pi/\alpha$ and $(k+1)\pi/\alpha \leq 4(w+1)$. Then $(-1)^k f(s, \alpha)$ is an analytic characteristic function on $V(k\pi/\alpha, (k+1)\pi/\alpha)$ provided $(-1)^w f(s)$, $(-1)^{k+w} m(z, \beta)$ are such for s on $V_{4w}$ and z on $V(k, k+1)$ respectively. This is detailed next. Say

$$(-1)^w f(s) = \int_R d(y)e^{sy}h_w(y) \text{ on } V_{4w} \text{ and}$$

$$(-1)^{k+w} m(z, \beta) = \int_R d(y)e^{zy}\omega_k(y, \beta) \text{ on } V(k, k+1),$$

with $h_w(y), \omega_k(y, \beta) > 0$ for all y. Assume that at each x with $k < x < k+1$, $e^{xy}\omega_k(y, \beta)$ is bounded for y real. Set $v_k(y, \alpha) := (1/4)(\pi/\alpha)^2\omega_k((\pi/\alpha)y, \pi^2/(4\alpha))$. Then

$$(-1)^k f(s,\alpha) = \int_R d(y)e^{sy}u_k(y, \alpha),$$

with $u_k(y, \alpha) := (v_k(\cdot, \alpha) * h_w)(y) := \int_R d(r)v_k(y-r, \alpha)h_w(r)$. Thus $u_k(y, \alpha) > 0$ for all real y.

**The Laplace representation of m(z, β).**

Assume $-\pi < \beta < \pi$. We will now determine the Laplace representation of m(z, β) for z on V(-1, 1) and that for z on V(k, k+1) when $k \neq 0, -1$.

**The Laplace representation of π/sin(πs).**

**Definition of $Q_k(z)$.** Let k be an integer. Set $Q_k(z) := z^{k+1}/(1+z)$ for $z \neq 0, -1$.

The partial fraction expansion of π/sin(πs) yields the following Laplace representation. If $0 < \mathrm{Re}(s) < 1$, then

$$\pi/\sin(\pi s) = \int_R d(y)e^{sy}/(1+e^y).$$

Note that for $x := \mathrm{Re}(s)$ one has $|e^{sy}| / (1+e^y)$ asymptotic with $e^{xy}$ as $y \to -\infty$ and with $e^{-(1-x)y}$ as $y \to \infty$. Thus $|e^{sy}| / (1+e^y)$ is integrable in y.

Now make the translation z = s+k. Then

$$\pi/\sin(\pi z) = (-1)^k \int_R d(y)e^{zy}Q_k(e^{-y}),$$

if $k < \mathrm{Re}(z) < k+1$. So π/sin (πs) is a meromorphic characteristic function on C.

When β is nonzero and z is not an integer, $m(z, \beta) = (1/(2i\beta))\sum_{\sigma=\pm1} \sigma\pi e^{i\sigma\beta z}/\sin(\pi z)$.

Assume $-\pi < \beta < \pi$ and $0 < \mathrm{Re}(z) < 1$.Then

$$\pi e^{i\beta z}/\sin(\pi z) = \int_R d(y)e^{zy}/(1+e^{y-i\beta}) .$$

This is a corollary of the principal established next, which concerns an instance where

$$e^{i\beta z}\int_R d(y)e^{zy}h(y) = \int_R d(y)e^{zy}h(y - i\beta),$$

or equivalently $\int_R d(y)e^{zy}h(y) = \int_R d(y)e^{z(y-i\beta)}h(y - i\beta)$.

**Translation principal** *Fix $\beta_0$, $\beta_1$ with $\beta_0 < \beta_1$. Assume h(γ) is analytic on an open set containing the closed horizontal strip $iV[\beta_0, \beta_1]$.*
**Definition of M(x).** Set $M(x) := \max\{|h(x+i\beta)| : \beta_0 \leq \beta \leq \beta_1\}$.
*Assume that $u_0 < u_1$ and $M(x) \sim O(\exp(-u_k x))$ as $x \to (-1)^{k+1} \cdot \infty$ on R. Say $u_0 < u_0' \leq u_1' < u_1$.*
(1) *There exists a $B'(u_0', u_1') \geq 0$ such that for all β, u with $\beta_0 \leq \beta \leq \beta_1$ and $u_0' \leq u \leq u_1'$, $\int_R d(y)e^{uy}|h(y + i\beta)| \leq B'(u_0', u_1')$. Say $z = u + it$ with $u_0 < u < u_1$. Assume $\beta_0 \leq \beta \leq \beta_1$. Then $|e^{zy}h(\gamma)|$ is integrable on the horizontal line through iβ.*
(2) *$I(\beta, z) := \int_R d(y)e^{z(y+i\beta)} h(y + i\beta)$ is constant in β: $I(\beta, z) = I(\beta_0, z)$.*

**Proof of (1).** $|h(y + i\beta)| \leq M(y)$. Also $e^{uy}M(y) \sim O(\exp((u-u_k)y))$ as $y \to (-1)^{k+1} \cdot \infty$ on R. $e^{uy} h(y + i\beta)$ is continuous in u, y, β on $D := (u_0, u_1)\times R\times[\beta_0, \beta_1]$ and therefore is bounded on compact subsets of D. Thus (1) holds.

**Proof of (2).** $e^{z\gamma}h(\gamma)$ is analytic in γ on S. Cauchy's theorem gives $\int d(\gamma)e^{z\gamma}h(\gamma) = 0$, with γ traversing once in the counterclockwise sense the rectangle R′ delineated next. Let $\beta$, $x_0$, $x_1$ with $\beta_0 \leq \beta \leq \beta_1$ and $x_0 < x_1$ be given. R′ has vertices $x_0 + i\beta_0$, $x_1 + i\beta_0$, $x_1 + i\beta$, $x_0 + i\beta$. Set $I(\beta, x_0, x_1, z) := \int d(y)e^{z(y+i\beta)} h(y + i\beta)$, with $x_0 \leq y \leq x_1$. Let $\varphi(x, \beta, z) := \int d(r)\, |e^{z(x+ir)} h(x + ir)|$, with $\beta_0 \leq r \leq \beta$. Then $|I(\beta, x_0, x_1, z) - I(\beta_0, x_0, x_1, z)| \leq \varphi(x_0, \beta, z) + \varphi(x_1, \beta, z)$. Now $z = u + it$ with $u_0 < u < u_1$ gives $\varphi(x, \beta, z) \leq (\beta - \beta_0)e^{|t|K} \cdot e^{ux}M(x)$, with $K := \max(|\beta_0|, |\beta_1|)$. So $\varphi(x, \beta, z) \to 0$ rapidly as $x \to \pm\infty$. Therefore (2) is valid.

Apply the translation $z = z'+k$ to obtain the following. Assume $-\pi < \beta < \pi$ and $k < Re(z) < k+1$. Then $\pi e^{i\beta z}/\sin(\pi z) = (-1)^k\int_R d(y)e^{zy}Q_k(e^{-(y-i\beta)})$.

This leads to the following lemma for the case $\beta \neq 0$, which we will use. The case $\beta = 0$ follows directly from the partial fraction representation of $\pi z/\sin(\pi z)$. It can also be obtained from the lemma with $\beta \neq 0$ by letting $\beta \to 0$.

**Definitions of $q(u, \beta)$, $l_k(y, \beta)$.** Set $q(u, \beta) := \sin(u\beta)/\beta$, if β is nonzero. Take $q(u, 0) := u$. Say β, y are real, β is not an odd multiple of π and k is an integer. Define $l_k(y, \beta) := e^{-ky}(e^{-y}q(k, \beta) + q(k+1, \beta))/(2(\cosh(y) + \cos(\beta)))$. Note that $l_k(y, -\beta) = l_k(y, \beta)$.

Fix β. $|e^{xy}l_0(y, \beta)| \sim O(e^{(|x|-1)|y|})$, when $|y|$ is large. If $|x| < 1$, then $e^{xy} l_0(y, \beta) \to 0$ rapidly as $y \to \pm\infty$ on $\mathbb{R}$. Thus $\int_R d(y)\, |e^{zy}l_0(y, \beta)|$ is finite, if $|Re(z)| < 1$. Say $k < x < k+1$ and y is real. Then $e^{xy}|l_k(y, \beta)|$ is $\sim O(e^{(x-(k+1))y})$ for $y \to \infty$ and is $\sim O(e^{(x-k)y})$ for $y \to -\infty$. So $k < Re(z) < k+1$ implies that $\int_R d(y)\, |e^{zy}l_k(y, \beta)|$ is finite.

**Lemma 1.2 Laplace representation of m(z, β).**
(1) *Assume $-\pi < \beta < \pi$ and $|Re(z)| < 1$. Then $m(z, \beta) = \int_R d(y)e^{zy}\, l_0(y, \beta)$, with $l_0(y, \beta) := q(1, \beta)/(2(\cosh(y) + \cos(\beta)))$ a probability density in y: $l_0(y, \beta) > 0$ and $\int_R d(y)\, l_0(y, \beta) = 1$.*
(2) *Assume*
*(*): $-\pi < \beta < \pi$ and k is an integer distinct from 0, -1.*
*Then, when $k < Re(z) < k+1$: $m(z, \beta) = (-1)^k \int_R d(y)e^{zy}l_k(y, \beta)$.*
(3) *Assume (*) and $\beta \neq 0$.*
*If w is an integer, $w\cdot(\pi / |\beta|) \leq k$ and $(k+1) \leq (w+1)\cdot(\pi / |\beta|)$, then $(-1)^w\cdot l_k(y, \beta) > 0$ for all real y.*

*If $k < w\cdot(\pi / |\beta|) < k+1$, then $l_k(y, \beta)$ changes sign exactly once for y on* **R**.
*If $k \neq 0, -1$, then $sign(k)\cdot l_k(y, 0) > 0$.*

**Corollary 1.1** Say *$\alpha > \pi/4$, $4w \leq k(\pi/\alpha)$ and $(k+1)(\pi/\alpha) \leq 4(w+1)$. $(-1)^k f(s, \alpha)$ is then an analytic characteristic function (ACF) on $V(k\pi/\alpha, (k+1)\pi/\alpha)$, provided $(-1)^w f(s)$ is an ACF for s on $V_{4w}$.*

**§2 The Laplace representation on $V_0'$ of f(s) := 1/(b(s)ζ(½ + s)) via that of 1/b(s).**

**Review** Part I, §3, Definition of $P_0(z)$ and Main unconditional theorem (1).

We will determine the unconditional representation when ½ < Re(s) < 4 of f(s) := 1/(b(s)ζ(½ + s)) as a two-sided Laplace transform: $f(s) = \int_R d(y)e^{sy}g(y)$.

First we proceed in a heuristic fashion and then we validate with rigor the promised results glimmering on intuition's horizon.

The only zeros z of the entire function n(s) := sin(πs/4)·2ξ(½ + s) with Re(z) > ½ are 4w with w ≥1. Those zeros are simple.

Say $w \geq 1$, $x_0 < x_1 \ldots < x_w$ and N(s) is a polynomial which has simple roots at $x_k$ for $1 \leq k \leq w$, but does not vanish elsewhere on the half-plane Re(s) > $x_0$. Then on $V(x_0, x_1)$: $1/N(s) = \int_R d(y)e^{sy}h(y)$ with $h(y) := -\sum_{1 \leq k \leq w} \exp(-x_k y)/N'(x_k)$ for y > 0.

In analogy with the polynomial case we guess that $g(y) = P_0(\pi e^{-2y})$ for y > 0.

Suppose that for s on $V_0'$, $1/b(s) = \int_R d(y)e^{sy}T_0(ie^{-2y})$. Now $1/\zeta(½ + s) = \sum_{n \geq 1}\mu(n) / n^{½ + s}$ when Re(s) > ½. Then

$$f(s) := 1/(b(s)\zeta(½ + s)) = \sum_{n \geq 1}(\mu(n)/n^{½}) \int_R d(y)e^{sy}T_0(i\cdot(ne^y)^{-2}).$$

If $\sum \int = \int \sum$ therein, then $g(y) = \sum_{n \geq 1}(\mu(n)/n^{½})T_0(i\cdot(ne^y)^{-2})$. So we expect that for 0 < v < 1:

$$P_0(\pi v) = \sum_{n \geq 1}(\mu(n)/n^{½})T_0(i\cdot v/(n^2)).$$

Now $P_0(\pi z) = -\sum_{k \geq 1}\sum_{n \geq 1}(\mu(n)/n^{½})(i\cdot z/(n^2))^{2k}/((\pi/4)a(½ + 4k))$. If $\sum_{k \geq 1}$, $\sum_{n \geq 1}$ can be interchanged therein, then $P_0(\pi z) = \sum_{n \geq 1} (\mu(n)/n^{½})(-\sum_{k \geq 1}(i\cdot z/(n^2))^{2k}/((\pi/4)a(½ + 4k)))$. So we guess that $T_0(z) = -(4/\pi)\sum_{k \geq 1} z^{2k}/a(½ + 4k)$, when z = ir with 0 < r < 1.

Stirling's formula (See G. Andrews, R. Askey, R. Roy [2]) shows that for $x \geq 2$, $a(x) \geq \pi^{1/2}\cdot x^{3/2}\cdot(x/(2\pi e))^{x/2}$. Thus $\sum_{k \geq 1} z^{2k} / a(½ + 4k)$ converges to an entire function.

**Definition of P(z).** Set $P(z) := (1/(\pi^{3/4}\Gamma(5/4)))\sum_{k \geq 1}(\pi z)^{2k}/((5/4)_{2k}(2k - ¼))$.

Our guess is that $T_0(z) = -P(z)$.

$$P(z) = (-4/(\pi^{3/4}\Gamma(5/4)))(-1 + (½)\sum_{\sigma = \pm 1} {}_2F_2(1, -¼ ; 5/4, ¾; \sigma\pi z)),$$

with ${}_2F_2$ a generalized hypergeometric function. See Eric W. Weisstein [15] for ${}_pF_q$.

See definition of $T_0(z, \beta)$ in §6.

**Lemma 2.1** *Say m is a nonnegative integer. Let $a_k$ be complex for k = m, m + 1, ... Assume $|a_k|^{1/k} \to 0$ as $k \to \infty$. Let E(z) be the entire function $\sum_{k \geq m} a_k\cdot z^k$. Say $q \geq 0$ and $Re(p) > 1 - mq$. Let Ω(z) be the entire function $\sum_{k \geq m} (a_k/\zeta(p + kq))\cdot z^k$. Then $\sum_{n \geq 1}(\mu(n)/(n^p))E(z/(n^q)) = \Omega(z)$.*

**Proof** Set $J(n, k, z) := (\mu(n)/(n^p))\cdot a_k\cdot(n^{-q} z)^k$. Then $|J(n, k, z)| \leq n^{-(Re(p)+mq)}|a_k|\cdot|z|^k$. Then $\sum_{n \geq 1}\sum_{k \geq m}|J(n, k, z)| \leq \zeta(Re(p) + mq) \sum_{k \geq m} |a_k|\cdot|z|^k < \infty$. So the interchange of summations $\sum_{n \geq 1}\sum_{k \geq m} J(n, k, z) = \sum_{k \geq m}\sum_{n \geq 1} J(n, k, z)$ holds.

**Corollary 2.1** *$P_0(\pi z) = \sum_{n \geq 1}(\mu(n)/(n^{1/2}))(-P(iz/(n^2)))$.*

See §6, Unconditional theorem 6.1, with β = ¼ and w = 0.

In validating the above heuristic argument we use the following lemmas. They are expressed using Mellin rather than Laplace transforms.

**Definitions of g(r, j, q), g(j, q)(r) and g(r, j, q, p).** Assume r > 0. Let j, q be real. Set $g(r, j, q) := r^q$, when $0 < r \leq 1$, and $g(r, j, q) := r^j$, when r > 1. Set g(j, q)(r) := g(r, j, q). Say p is real. Define $g(r, j, q, p) := \sum_{n \geq r} n^{-p}g(r/n, j, q) \leq \infty$.

g(r, j, q) is continuous. g(r, j, q) > 0. Also $g(r, j, q, p) = r^j \cdot\sum_{1 \leq n < r} n^{-(p+j)} + r^q \cdot\sum_{n \geq r} n^{-(p+q)}$. Note that g(r, j, q, p) is finite, if p + q > 1. When $0 < r \leq 1$, $g(r, j, q, p) = \zeta(p + q)\cdot r^q$.

**Lemma 2.2** *Assume p + q > 1, j ≠ 1 – p and r > 1. Then:*
*$g(r, j, q, p) \leq \alpha(p + j)\cdot r^j + \beta(j, q, p)\cdot r^{1-p}$, with α(u) := 1 + (1/(u - 1)), for u ≠ 1, and β(j, q, p) := 1 + (1/(p + q - 1)) - (1/(p + j - 1)).*

**Proof** The lemma is a consequence of the Claims 1, 2 established next.

Claim 1 *Assume $u > 1$ and $r \geq 1$. Then $\sum_{n \geq r} n^{-u} \leq \alpha(u)\cdot r^{-(u-1)}$.*

Proof of Claim 1. $\sum_{n \geq m+1} n^{-u} \leq \int_{x \geq m} d(x)x^{-u}$, with m the least integer $\geq r$.

Claim 2 *Assume $u > 0$, $u \neq 1$ and $r > 1$. Then*

$$\sum_{1 \leq n < r} n^{-u} \leq 1 + (1/(u-1))(1 - r^{-(u-1)}).$$

Proof of Claim 2. $\sum_{2 \leq n < r} n^{-u} \leq \int_{1 \leq x < r} d(x)x^{-u}$.

**Corollary 2.2** Assume $p + q > 1$ and $j \neq 1 - p$. Then $g(r, j, q, p) \leq K(j, q, p)g(r, j', q)$, with $j' = \max\{1 - p, j\}$ and $K(j, q, p) = \max\{\zeta(p + q), |\alpha(p + j)| + |\beta(j, q, p)|\}$.

**Lemma 2.3**
(1) *If $j < x < q$, then $\int_{v>0} (dv)v^{x-1}g(1/v, j, q)$ is finite.*
*Assume $p + q > 1$, $j \neq 1 - p$ and $\max\{1 - p, j\} < x < q$.*
(2) *$\int_{v>0} (dv)v^{x-1}g(1/v, j, q, p)$ is finite.*
(3) *$\zeta(p + x)\cdot\int_{v>0} (dv)v^{x-1}g(1/v, j, q) = \int_{v>0} (dv)v^{x-1}g(1/v, j, q, p)$.*

**Definitions of $q_0(T)$ and $j_0(T)$.** Say $r > 0$. Let $T(r)$ be a complex-valued measurable function of r. Assume that there is a real q such that for some nonnegative $K(q)$: $|T(r)| \leq Kr^q$, when $0 < r \leq 1$. Let $q_0(T)$ be the least upper bound of all such q. Then $q_0(T) \leq \infty$. Suppose there is a real j such that for some $K_1(j) \geq 0$: $|T(r)| \leq K_1(j)r^j$, when $r > 1$. Let $j_0(T)$ the greatest lowest bound of all such j. Then $j_0(T) \geq -\infty$.

Say $q < q_0(T)$ and $j > j_0(T)$.

**Lemma 2.4** *Assume that there exist real j, q and a $K(j, q) \geq 0$ such that $|T(r)| \leq K(j, q)g(r, j, q)$.*
(1) *If $j_0(T) < Re(s) < q_0(T)$, then $\int_{v>0} (dv)|v^{s-1}T(1/v)|$ is finite.*
(2) *Let p be complex. Assume $Re(p) + q_0(T) > 1$.*
(i) *$\sum_{n \geq 1} n^{-Re(p)}|T(r/n)|$ is finite.*
**Definitions of $\theta_{T,p}(r)$ and $\omega_{T,p}(r)$.** Let

$$\theta_{T,p}(r) := \sum_{n \geq 1} (n^{-p})T(r/n) \text{ and } \omega_{T,p}(r) := \sum_{n \geq 1} (\mu(n)/(n^p))T(r/n).$$

(ii) *Say $T(r)$ is continuous for $r > 0$. Then so are $\theta_{T,p}(r)$ and $\omega_{T,p}(r)$.*
*Let $\varphi = \theta_{T,p}$ or $\varphi = \omega_{T,p}(r)$.*
(ii)′ *Say p is real and $T(r) > 0$, for all $r > 0$. $\theta_{T,p}(r) > 0$.*
*Set $p_1 = Re(p)$.*

(iii) *$j_0(\varphi) \leq max\{1 - p_1, j_0(T)\}$ and $q_0(\varphi) \geq q_0(T)$.*
*Assume $max\{1 - p_1, j_0(T)\} < Re(s) < q_0(T)$.*
(iii)′ *$\int_{v>0} (dv)|v^{s-1}\varphi(1/v)|$ is finite.*
(iv)

$$\zeta(p + s)\int_{v>0} (dv)v^{s-1}T(1/v) = \int_{v>0} (dv)v^{s-1}\theta_{T,p}(1/v).$$

(iv)′

$$(1/\zeta(p + s))\int_{v>0} (dv)v^{s-1}T(1/v) = \int_{v>0} (dv)v^{s-1}\omega_{T,p}(1/v).$$

Note that $\int_{y<0} (dy)|e^{sy}P_0(\pi e^{-2y})|$ is finite when Re(s) > ½ , provided that for small $\varepsilon > 0$, $|P_0(r)| \sim O(r^{1/4+\varepsilon})$ for $r > \pi$.

**Corollary 2.3** *Assume $\varepsilon > 0$ and $|P(ir)| = O(r^{1/4+\varepsilon})$ for $r > 1$. Then:*
*$|P_0(\pi r)| = O(r^{1/4+\varepsilon})$, for $r > 1$.*

## §3 Strategy for determining the Laplace representation on $V_0'$ of 1/b(s).

We now outline the strategy for determining the two-sided Laplace transform representation of 1/b(s).

**Divergence of the formal partial fraction expansion of 1/b(s).**

We do not use the formal partial fraction expansion of 1/b(s) to obtain its Laplace transform representation. That expansion is D(s) + C(s) with: $D(s) := \sum_{w \leq -1} (1/b'(4w))(1/(s - 4w))$ and $C(s) := j(s) + (1/b'(0))(1/s) + (1/b'(½))(1/(s - ½))$, where $j(s) := \sum_{w \geq 1} (1/b'(4w))(1/(s - 4w))$.

Say s ≠ 4w for any integer w ≥ 0 and s ≠ ½. Then the series for j(s) converges and C(s) is analytic. That is a result of the following. Given d > 0, let K(d) be the set of s with |s - 4w| ≥ d for all integers w ≥ 1. The series for j(s) converges uniformly in s on K(d). The reason is that for w ≥ 1 and s in K(d): $|(1/b'(4w))(1/(s - 4w))| \leq \pi^{2w}/((2w)!\cdot d)$.

However the formal expansion D(s) diverges for all s, since 1/b′(4w) → ∞ ultrarapidly as the integer w → -∞. The latter is seen as follows.

$b'(4w) = (-1)^w\pi^{3/4 - 2w}\Gamma(5/4 + 2w)(2w - 1/4)$. Say w ≤ -1. Now Γ(s + 1) = s Γ(s). Thus when k is a nonnegative integer: $1/\Gamma(\varepsilon - (k + 1)) = (-1)^{k+1}((1 - \varepsilon)/\Gamma(\varepsilon))(2 - \varepsilon)_k$. Set ε = ¼ and k = -2(w + 1). Then $1/b'(4w)) = (-1)^w\pi^{-(2 + 3/4)}((1 - \varepsilon)/\Gamma(\varepsilon))E(k)$ with $E(k) := \pi^{-k}\cdot(2 - \varepsilon)_k/(k + 2 + 1/4)$. Now ε < 1 implies $(2 - \varepsilon)_k \geq k!$. Also $k! \geq e\cdot(k/e)^k$.

**Definitions of $n_0(s)$ and $f_0(s)$.** Set

$$n_0(s) := \sin(\pi s/4)\mathrm{l}(½ + s) \text{ and } f_0(s) := 1/n_0(s).$$

$b(s) := (s - ½)n_0(s)$. The formal derivation of the Laplace representation on $V_0'$ of $1/b(s)$ from

$$f_0(s) = \int_R d(y)e^{sy}H_0(e^{-2y})$$

via a convolution gives

$$1/b(s) = \int_R d(y)e^{sy}T_0(ie^{-2y}), \text{ with } T_0(iJ) = J^{¼} \int_{0 < j \leq J} d(j)j^{-5/4}(½)H_0(j), \text{ for } J > 0.$$

See §6, Corollary 6.1, with $\beta = ¼$ and $w = 0$.

The next lemma will reduce aspects of the study of the representation of $1/b(s)$ to that of $f_0(s)$.

Say $j > 0$. Let $h(j)$ be a measurable complex–valued function. Suppose that for some positive $\varepsilon$, $|h(j)| \leq K\cdot j^q$ for j with $0 < j \leq \varepsilon$ . Say $|h(j)|$ is bounded on any interval $[\varepsilon, J]$. Assume $\mathrm{Re}(\alpha) < q$. Say $J > 0$. $\int_{0 < j \leq J} d(j)|j^{-(1+\alpha)}h(j)|$ is finite.

**Definition of $h^{<\alpha>}$.** Set $h^{<\alpha>}(J) := J^{\alpha}\int_{0 < j \leq J} (d(j)/j)j^{-\alpha}h(j)$.

**Lemma 3.1** *Assume $\mathrm{Re}(\alpha) < q$, j is real and $j \neq \alpha$. Take $m = \max\{\mathrm{Re}(\alpha), j\}$. There is a $K(\alpha, q, j) \geq 0$ such that $0 < (g(j, q))^{<\alpha>}(J) \leq K(\alpha, q, j)\cdot g(m, q)(J)$, for all $J > 0$.*

**Proof** If $0 < J \leq 1$, then $(g(j, q))^{<\alpha>}(J) = (1/(2(q - \alpha)))\cdot J^q$. If $J > 1$, then $(g(j, q))^{<\alpha>}(J) = (1/(2(j - \alpha)))\cdot(((j - q)/(q - \alpha))\cdot J^{\alpha} + J^j)$ .

**Lemma 3.2** *Assume $q > \mathrm{Re}(\alpha)$ and $|h(r)| \leq Kg(j, q)(r)$, for $r > 0$.*
(1) *Take $m = \max\{\mathrm{Re}(\alpha), j\}$. There is a $K' \geq 0$ such that $|h^{<\alpha>}(J)| \leq K'\cdot g(m, q)(J)$, for all $J > 0$. $j_0(h^{<\alpha>}) \leq \max\{\mathrm{Re}(\alpha), j_0(h)\}$. $q_0(h^{<\alpha>}) \geq q_0(h)$.*
(2) *Assume $q_0(h) > \mathrm{Re}(\alpha)$. Say $\max\{\mathrm{Re}(\alpha), j_0(h)\} < \mathrm{Re}(s) < q_0(h)$. Then $\int_R d(y)|e^{sy}\theta(e^{-y})| < \infty$, for $\theta = h, h^{<\alpha>}$.*

$$(1/(s - \alpha))\int_R d(y)e^{sy}h(e^{-y}) = \int_R d(y)e^{sy}h^{<\alpha>}(e^{-y}).$$

(3) *Suppose that $h(z)$ is an entire function of order q at 0. Then, for $z > 0$,*

$$h^{<\alpha>}(z) = \sum_{n \geq q} (h^{(n)}(0)/(n!))(1/(n - \alpha))\cdot z^n.$$

*Thus $h^{<\alpha>}$ extends to an entire function.*
(4) *If $h(j) > 0$ for $j > 0$, then $h^{<\alpha>}(j) > 0$ for $j > 0$.*

**Corollary 3.1** *Suppose that $|h(r)| \le Kg(j, q)(r)$, for $r > 0$. Set $m = \max\{1 - Re(p), Re(\alpha)\}$. Assume $q_0(h) > m$. Let $T = h^{<\alpha>}$. $|\omega_{T,p}| \le Kg(\max\{m, j\}, q)$. Then $j_0(\omega_{T,p}) \le \max\{m, j_0(h)\}$ and $q_0(\omega_{T,p}) \ge q_0(h)$.*

**Proof** Apply (1), (2) of the previous Lemma 3.2 and (2) (i), (ii) of Lemma 2.4 of §2.

Next Theorem 3.1 of Part I, §3, is proven.

Say w is an integer. f(s) is analytic on $V_{4w}'$. Say x and t are real.

**Integrability of |f(s)| on vertical lines outside of the critical strip.**

**Theorem 3.1** *|f(s)| is both integrable and square integrable on each vertical line $x + i\mathbf{R}$ with x real, $|x| \ge ½$ and x not a nonzero multiple of four.*

**Proof** f(s) is odd. So assume $x \ge ½$ and $x \ne 4w$ for $w \ge 1$. We apply the asymptotic expansion of $|b(x + it)|$ of Part I, §4, (4.1). $1 / |b(x + it)|$ is integrable (respectively: square integrable) in t over $\mathbb{R}$, if $x > -3/2$ (respectively: $x > -5/2$) and $x \ne 4w$, $w \ge 0$.

If $v \ge \sigma > 1$, then $|\zeta(v + it)| \ge \zeta(2\sigma)/\zeta(\sigma)$. That has the following consequence. Say $½ < u < x < x_1$. When $|t| \ge T > 0$,

$$|f(x + it)| \le K(u, x_1)\cdot|t|^{-(7/4 + x/2)} \le K(u, x_1)\cdot|t|^{-(7/4 + u/2)}.$$

So the theorem holds for $x > ½$.

The case $x = ½$ and the case $x > ½$ are consequences of the theorem stated next and proven in Apostle [3], ch. 13, p 287.

**Bound for $1 / |\zeta(z)|$.** There exists an $M > 0$ such that for any z with $Re(z) \ge 1$ and $|t| \ge e$, where $t := Im(z)$, one has $1 / |\zeta(z)| < M\cdot(\log(|t|))^7$.

Then $|f(x + it)| \le K(u, x_1)\cdot(\log(|t|))^7\cdot|t|^{-(7/4 + u/2)}$, when $½ \le u < x < x_1$ and $|t| \ge e$.

Say β is complex.

**Definitions of N(z, β) and F(z, β).**

$$N(z, \beta) := (1/\pi)\sin(\pi z/2)\cdot\Gamma(1 + \beta + z).\ F(z, \beta) := 1/N(z, \beta).$$

F(z, β) is analytic in z except for the following simple poles. When β is not an integer, each even integer is a pole. When β is an integer, the poles are the even $n - \beta$ arising from the integers $n \geq 0$. When z is not a pole, F(z, β) is entire in β.

$n_0(s) = 2(\pi^{3/4 - s/2})N(s/2\ , 1/4)$. Hence

$$f_0(s) = \tfrac{1}{2}\cdot(\pi^{-3/4 + s/2})F(s/2\ , 1/4).$$

**Integrability of |N(z, ¼)| on vertical lines.** |N(z, ¼)| is integrable (respectively: square integrable) on each vertical line $x + i\mathbf{R}$ with $x > 1/4$ (respectively: $x > -1/4$) and $x \neq 2w$. This is a consequence of the asymptotic expansion of N(z, ¼) for $z = x + it$, with $x_0 < x < x_1$, t real and |t| large.

$$|N(x + it, 1/4)| \sim (2\pi)^{1/2}\cdot|t|^{3/4 + x/2}(1 + \varepsilon(x,t)/|t|),$$

with $|\varepsilon(x,t)| < K(x_0, x_1) < \infty$.

**Determination of the Laplace representations of $f_0(s)$ and 1/b(s) on $V_0'$ from that of F(z, ¼) on V(¼, 2).** Say H(r) is continuous for $r > 0$. Assume: $H(r) = O(r^2)$, for all small r; and for arbitrarily small positive ε, $H(r) = O(r^{1/4 + \varepsilon})$, for all large r. Suppose that for z on V(¼, 2):

$$F(z, 1/4) = \int_R d(y)e^{zy}H(e^{-y}).$$

Then on $V_0'$,

$$f_0(s) = \int_R d(y)e^{sy}H_0(e^{-2y}),\ \text{with } H_0(v) := \pi^{-3/4}H(\pi v).$$

Apply Lemma 3.2. On $V_0'$,

$$1/b(s) = \int_R d(y)e^{sy}T_0(ie^{-2y}),\ \text{with } T_0(iJ) = \tfrac{1}{2}H_0^{<1/4>}(J).$$

**Translation relation for F(z, β).** Say w is a nonnegative integer.

$$(-1)^w F(z + 2w, \beta) = (1/(1 + \beta + z)_{2w})F(z, \beta).$$

**Definition of j(u, m).** Set $j(u, m) := 1/(\zeta(u + 2m)(1 + (\tfrac{1}{2})u)_m)$.

$j(u, m)$ is a meromorphic function on C. Say $m \geq 2$. Assume $Re(u) > -3$. Then $Re(u + 2m) > 1$. Say $Re(u) > -2$. Then $j(u, m)$ is analytic in u.

**Translation relation for f(z).** Say $w \geq 1$. Set $m = 2w$. The translation relation for $F(z, ¼)$ gives

$$(-1)^w f(S + 2m) = \pi^m(1/(S + 2m - ½))j(½ + S, m)f_0(S).$$

**The Mellin transform representation of j(u, m).**

**The Mellin transform representation of $1/(z)_m$.** Say m is an integer $\geq 1$. The partial fraction expansion of $1/(z)_m$ implies that when $Re(z) > 0$:

$$1/(z)_m = \int_{0 < v \leq 1} (dv)v^{z-1}(1/((m-1)!))(1 - v)^{m-1}.$$

That is a special case of the representation of the Beta function:

$$\Gamma(p)\Gamma(q)/\Gamma(p + q) = \int_{0 < v \leq 1} (dv)v^{p-1}(1 - v)^{q-1},$$

when $Re(p) > 0$ and $Re(q) > 0$. Thus, when $Re(u) > -2$, one has:

$$1/(1 + (½)u)_m = \int_{0 < v \leq 1} (dv)v^u(2v/((m-1)!))(1 - v^2)^{m-1}.$$

Say $m \geq 2$ and $0 < v \leq 1$.

**Definition of E(v, m).** Set $E(v, m) := (2/((m-1)!))\cdot v^2\sum_{1 \leq k < 1/v} \mu(k)(k^{-2} - v^2)^{m-1}$.

**Lemma 3.3 Mellin transform representation of j(u, m).** *Say $m \geq 2$.*
(1) *If $0 < v < 1$, then*

$$0 < ((2 - \pi^2/6)/((m-1)!))\cdot v^2(1 - v^2)^{m-1} < E(v, m) < (\pi^2/6)(2/((m-1)!))\cdot v^2.$$

(2) *On the half-plane of u with $Re(u) > -2$: $j(u, m) = \int_{0 < v \leq 1} (dv)v^{u-1}E(v, m)$.*

**Proof of (1).** $|\mu(k)(k^{-2} - v^2)^{m-1}| < k^{-2(m-1)}$. Thus $\sum_{1 \leq k < 1/v} |\mu(k)(k^{-2} - v^2)^{m-1}| < \zeta(2(m - 1)) \leq \pi^2/6$. Also $E(v, m) = (2/((m-1)!))\cdot v^2(1 - v^2)^{m-1}(1 + e(k, v))$, with $e(k, v) := \sum_{2 \leq k < 1/v} \mu(k)((k^{-2} - v^2)/(1 - v^2))^{m-1}$. Now $0 < ((k^{-2} - v^2)/(1 - v^2))^{m-1} < k^{-2(m-1)}$. Thus $1 + e(k, v) > 2 - \zeta(2(m - 1)) \geq 2 - \pi^2/6$.

**Proof of (2).** We apply the Lemma 2.4 (2), (iv)′ of §2 as follows. Take $T(r) := (2/((m-1)!))\cdot r^{-2}(1 - r^{-2})^{m-1}$, for $r > 1$. Set $T(r) := 0$, when $0 < r \leq 1$. Then $q_0(T) = \infty$ and $j_0(T) = -2$. Let $p = 2m$. So $j(u, m) = \int_{0 < v \leq 1} (dv)v^{u-1}E(v, m)$, when $Re(u) > -2$.

**The Mellin transform representations of 1/(S + 2m - ½), j(½ + S, m) and $f_0$(S).**

Assume that Re(S) > -2m + ½. Then $1/(S + 2m - ½) = \int_{0 < v \leq 1} (dv)v^{S-1}g_1(v, m)$, with $g_1(v, m) := v^{2m - ½}$.

Assume that Re(S) > -(2 + ½). Then $j(½ + S, m) = \int_{0 < v \leq 1} (dv)v^{S-1}g_2(v, m)$, with $g_2(v, m) := v^{½}E(v, m)$.

Assume that for ½ < Re(S) < 4, $f_0(S) = \int_{v > 0} d(v)v^{S-1}H_0(v^{-2})$, with $q_0(H_0) \geq 2$ and $j_0(H_0) \leq ¼$.

The last integral converges absolutely. (See the case β = ¼ of Lemma 5.3 (4) and of Corollary 6.1 (1).)

**The Mellin transform representation of f(s), for s on $V_0'$ + 4w.**

$$(-1)^w f(S + 2m) = \int_{v > 0} d(v)v^{S-1}G(v, m),$$

where

$$G(v, m) := \pi^m \iint_D d(v_1)d(v_2)v_1^{-1}g_1(v, m)v_2^{-1}g_2(v, m)H_0((v_1v_2/v)^2), \text{ with } D := (0,1]^2.$$

The last two integrals converge absolutely. G(v, m) is continuous in the v > 0.

Then for s on $V_0'$ + 4w,

$$f(s) = \int_{v > 0} d(v)v^{s-1}v^{-2m}G(v, m).$$

**Positivity** G(v, m) is positive for all v > 0, provided that holds for H(v). $v^{-2m}G(v, m)$ is then positive. (See the case β = ¼ of the Unconditional theorem 6.1 (2) and of Lemma 6.2 (2).)

**Determination of the Laplace representation of f(s) on $V_{4w}'$.** In order to determine the Laplace representation of f(s) on $V_{4w}'$, for w ≥ 0, we will determine that of F(z, β) relative to z on $V_{4w}$ for certain β including β = ¼. (See Unconditional theorem 6.1 (1), (2).) The analytic dependence on β of expressions in certain relations which arise and the extensions thereby of some of them play an essential role, even in obtaining the desired results for the case β = ¼.

**Determination of the Laplace representation of F(z, β) relative to z on $V_{4w}$.**

**Integrability of |F(z, β)| on vertical lines.** Assume β is complex. Let x, t be real. $|F(x + it, β)| \sim 0(|t|^{-(½ + Re(β) + x)})$, when |t| is large. |F(x + it, β)| is square

integrable in t on **R** for arbitrarily small positive x precisely when $\mathrm{Re}(\beta) \geq 0$. If $\mathrm{Re}(\beta) \geq 0$, then for each $x > 0$ which is not a pole of $F(z, \beta)$, $|F(x + it, \beta)|$ is square integrable in t.

**Divergence of the formal partial fraction expansion of F(z, β).** We seek to determine the Laplace representation of $F(z, \beta)$ relative to z on the strip $V(2)$: $0 < \mathrm{Re}(z) < 2$. We do not employ the formal partial fraction expansion of $F(z, \beta)$. Say β is not an integer. That expansion $\sum_{-\infty < w < \infty}(1/N_z(2w, \beta))(1/(z - 2w))$ has a divergent component as is shown next. $N_z(2w, \beta) = (-1)^w \cdot ½\Gamma(1 + \beta + 2w)$. Thus $1/N_z(2w, \beta) = (-1)^{w+1}2(\sin(\pi\beta)/\pi)\Gamma(-\beta + 2k)$, with $k = -w$. Now $\sin(\pi\beta) \neq 0$. Also $|\Gamma(-\beta + 2k)| \to \infty$ as $k \to \infty$. The latter is an instance of the following. Say $\delta \geq 0$ and $r \geq 1$. Then for any $z = x + iy$ in the semi-strip $x \geq r$, $-\delta \leq y \leq \delta$: $|\Gamma(x + iy)| \geq \theta(\delta)\Gamma(r)$, with $\theta(\delta) := (\pi\delta/\sinh(\pi\delta))^{½}$, for $\delta > 0$, and $\theta(0) = 1$. Thus $|\Gamma(z)|$ diverges to infinity with uniform ultrapidity as z tends to infinity in the semi-strip.

We will determine the representation of $F(z, \beta)$ relative to z on $V(2)$ via that of the following function $F(u, \beta, 1)$, with $-(1 + \beta) < \mathrm{Re}(u) < 1 - \beta$.

**Definitions of N(u, β, 1) and F(u, β, 1).** Take $u = 1 - (\beta + z)$.

$$N(u, \beta, 1) := (1/\pi)\cos((\pi/2)(u + \beta))\cdot\Gamma(2 - u).\ F(u, \beta, 1) := 1/ N(u, \beta, 1).$$

Then $N(z, \beta) = N(u, \beta, 1)$ and $F(z, \beta) = F(u, \beta, 1)$.

$$F(u, \beta, 1) = (\sin(\pi u)/\cos((\pi/2)(u + \beta)))\cdot(1/(1 - u))\Gamma(u),$$

since $1/(\Gamma(u)\Gamma(1 - u)) = (1/\pi)\sin(\pi u)$. So we will use the Laplace representation of $\Gamma(u)$ for u with $\mathrm{Re}(u) > 0$.

**Splitting F(u, β, 1).**

Next we split $F(u, \beta, 1)$ as a sum of two functions. For the ranges of $\mathrm{Re}(u)$, $\mathrm{Re}(\beta)$ to be considered, one of those functions decays with exponential rapidity and the other with algebraic rapidity as $\mathrm{Im}(u) \to \pm\infty$.

**Trigonometric identity.** We use

$$\sin(\pi u)/\cos((\pi/2)(u + \beta)) = \sin(\pi\beta)/\cos((\pi/2)(u + \beta)) + 2\sin((\pi/2)(u - \beta)),$$

which is obtained as follows. $2\sin(\theta - \omega)\cos(\theta + \omega) = \sin(2\omega) - \sin(2\theta)$. Hence

$$\sin(2\theta)/\cos(\theta + \omega) = \sin(2\omega)/\cos(\theta + \omega) + 2\sin(\theta - \omega).$$

**Definitions of E(u, β, 1) and E(u, β, 2).**

$$E(u, \beta, 1) := (1/(1 - u))\cdot((\pi/2)/\cos((\pi/2)(u + \beta)))\cdot\Gamma(u).$$

$$E(u, \beta, 2) := 2\sin((\pi/2)(u - \beta))\cdot(1/(1 - u))\Gamma(u).$$

**The splitting of F(u, β, 1) via E(u, β, 1) and E(u, β, 2).**

$$F(u, \beta, 1) = (2/\pi)\sin(\pi\beta)E(u, \beta, 1) + E(u, \beta, 2).$$

We will use the splitting of F(u, β, 1) to obtain its Mellin transform representation of Lemma 3.8 from that of E(u, β, 2) in Corollary 3.2 and that of E(u, β, 1) in Lemma 3.7.

**Trigonometric asymptotics.** The asymptotic behavior of $|2\sin((\pi/2)(u - \beta))|$ and of $|((\pi/2)/\cos((\pi/2)(u + \beta)))|$, for $|\mathrm{Im}(u)|$ large, is determined by the next observations. Henceforth take

$$\beta = \beta_1 + i\beta_2, \text{ with } \beta_k \text{ real.}$$

Let x, t be real. Take θ(z) to be either of the functions sin(z) or cos(z). $|2\theta(x + it)| \sim e^{-|t|}(1 + B(x, t))$, with $|B| \leq e^{-2|t|}$. Let t, r be real. $||t + r| - |t|| \leq |r|$. Say t is nonzero. Let $\sigma := \mathrm{sign}(t)$. If $|t| \geq -\sigma r$, then $|t + r| - |t| = \sigma r$.

**The asymptotic behavior of |E(u, β, 1)| and of |E(u, β, 2)|.** Assume $x_0 < x < x_1$, $|\beta| \leq B$, $|t + \beta_2| > \delta > 0$ and $|t| > T > 0$ with T large. Say u = x + it. Apply the previous trigonometric asymptotics and Stirling's approximation to Γ(u). Then the following approximations hold uniformly. $|((\pi/2)/\cos((\pi/2)(u + \beta)))\cdot\Gamma(u)| \sim O(e^{-\pi|t|}\cdot|t|^{-(½- x)})$. Hence $|E(u, \beta, 1)| \sim O(e^{-\pi|t|}\cdot|t|^{-(3/2- x)})$. Also $|(1/(1 - u))\Gamma(u)| \sim O(e^{-\pi|t|/2}\cdot|t|^{-(3/2- x)})$. So $|E(u, \beta, 2)| \sim O(t|^{-(3/2- x)})$.

**Integrability.**

$\int_R|\theta(x + it)|^p\, dt$ is finite for any p > 0 in the following cases.
(1) $\theta(u) := (1/(1 - u))\Gamma(u)$ and x is not an integer ≤ 1.
(2) Assume (*): x is not an integer ≤ 0 and x ≠ - β + n with n an odd integer. $\theta(u) := ((\pi/2)/\cos((\pi/2)(u + \beta)))\cdot\Gamma(u)$.
(3) Say (*) and x ≠ 1 hold. θ(u) := E(u, β, 1). Assume β is not an integer. Fix β. $\theta(u) := (2/\pi)\sin(\pi\beta)E(u, \beta, 1)$.

Say β is not an integer. |E(u, β, 2)| is square integrable (respectively integrable) on any vertical line x + i**R** with x < 1 (respectively x < ½) and x not a negative

integer or zero.

**The Mellin transform representation of E(u, β, 2).**

**Definition of J(z).** Set $J(z) := (1- e^{-z})/z$, for z nonzero. Take $J(0) := 1$.

J(z) is an entire function. Assume $z \neq 0$. If $Re(z) \geq 0$, then $|J(z)| \leq 2 / |z|$. Also: $|J(z)| \leq (e^{|z|} - 1) / |z| \leq e^{|z|}$.

**Lemma 3.4** *Assume $0 < Re(u) < 1$.*
(1) $\int_{v > 0} (dv)\ |v^{u-1}J(v)|$ *is finite.*
(2) $(1/(1 - u))\Gamma(u) = \int_{v > 0} (dv) v^{u-1}J(v)$.

**Lemma 3.5** *Assume $0 < Re(u) < 1$ and $-\pi/2 \leq \varphi \leq \pi/2$.*
$e^{i\varphi u}(1/(1 - u))\Gamma(u) = \int_{\mathbf{R}} d(y)e^{uy}J(e^{y} e^{-i\varphi})$.

**Proof** We show that $e^{i\varphi u} \int_{\mathbf{R}} d(y)e^{uy}J(e^{y}) = \int_{\mathbf{R}} d(y)e^{uy}J(e^{y - i\varphi})$. The translation principle of **§1** applies as follows.
Say $v > 0$. Set $M(v) := \max\{|J(ve^{i\theta})| : -\pi/2 \leq \theta \leq \pi/2\}$. $Re(ve^{i\theta}) \geq 0$. $M(v) \leq 2/v$. Also $M(v) \leq e^{v}$. So $M(e^{j}) \leq e$, for $j \leq 0$.

**Corollary 3.2** *Assume $0 < Re(u) < 1$.*
(1) *Say $-\pi/2 \leq \theta \leq \pi/2$ and ω is real.* $\int_{v > 0} d(v)|v^{u-1}Im(e^{-i\omega}J(ve^{-i\theta}))|$ *is finite. Also* $\sin(\theta u - \omega)(1/(1 - u))\Gamma(u) = \int_{v > 0} d(v)v^{u-1}Im(e^{-i\omega}J(ve^{-i\theta}))$.
(2) *Say ω is complex.* $\int_{v > 0} d(v)|v^{u-1}(1/v)(\cos(\omega) - \cos(v - \omega))|$ *is finite. Also:*
**The Mellin transform representation of E(u, β, 2).**

$$½E(u, (2/\pi)\omega, 2) = \sin((\pi/2)u - \omega)(1/(1 - u))\Gamma(u)$$

$$= \int_{v > 0} d(v)v^{u-1}(1/v)(\cos(\omega) - \cos(v - \omega)).$$

**The Mellin transform representation of E(u, β, 1).**

Say $-1 < Re(u) < 1$. Then $(\pi/2)/\cos((\pi/2)u) = \int_{v > 0} d(v)v^{u-1}(1/(v + 1/v))$. Now assume $-(1 + \beta_1) < Re(u) < 1 - \beta_1$. Then

$$(\pi/2)/\cos((\pi/2)(u + \beta)) = \int_{v > 0} d(v)v^{u-1}v^{\beta}(1/(v + 1/v)).$$

Suppose that $z > 0$ or z is nonreal. $z^{\beta}$ is obtained below using the principal branch $-\pi < \arg(z) < \pi$.

**Definition and properties of W(z, β).**

**Claim 3.1** *Assume either of (1), (2) below.*
(1) $\beta_1 < 1$ *and* $Re(z) > 0$.
(2) $-1 < \beta_1 < 1$ *and* $Re(z) \geq 0$.
*Then the integral*

$$W(z, \beta) := \int_{j>0} d(j) j^{-\beta}(1/(1+j^2))e^{-zj}$$

*converges absolutely.*
*Say (1) holds. Then* $W(z, \beta)$ *is analytic in each of* $z, \beta$.

**Translation relation for W(z, β).**

$$W(z, \beta) = \Gamma(1-\beta)z^{\beta-1} - W(z, \beta-2).$$

*Assume (2) holds. Then* $W(z, \beta)$ *is analytic in* $\beta$ *and continuous in* $z$. *Also* $W(0, \beta) = (\pi/2)/\cos((\pi/2)\beta)$.

**Proof** Suppose $j \geq 1$. In case (1), say r is real and $\delta > 0$. Consider $\beta$, z with $\beta_1 \geq r$ and $Re(z) \geq \delta$. Then $|j^{-\beta}(1/(1+j^2))e^{-zj}| \leq j^{-r-2}e^{-\delta j}$. In case (2), take r with $r > -1$. If $\beta_1 \geq r$, then $|j^{-\beta}(1/(1+j^2))e^{-zj}| \leq j^{-r-2}$.
Say $0 < j < 1$ and $r < 1$. Consider $\beta$ with $\beta_1 \leq r$. Then $|j^{-\beta}(1/(1+j^2))e^{-zj}| \leq j^{-r}$.
If (1) holds then the translation relation for $W(z, \beta)$ is valid, since $1/(1+j^2) = 1 - (j^2/(1+j^2))$.

**Lemma 3.6** *Assume* $\beta_1 < 1$ *and* $\max\{0, -(1+\beta_1)\} < Re(u) < 1 - \beta_1$.
(1) $\int_{v>0} d(v)|v^{u-1}W(v, \beta)|$ *is finite.*
(2) **Mellin transform representation of ((π/2)/cos((π/2)(u + β)))·Γ(u).**

$$((\pi/2)/\cos((\pi/2)(u+\beta)))\cdot\Gamma(u) = \int_{v>0} d(v)v^{u-1}W(v, \beta).$$

**Proof** Say $v > 0$. Let $e_0(v, \beta) := v^\beta(1/(v + 1/v))$. The convolution $\int_{j>0} d(j)j^{-1}e_0(j/v, \beta)e^{-j}$ equals $W(v, \beta)$, since $e_0(r, \beta) = e_0(1/r, -\beta)$.

**Definition and properties of $B_0(z, \beta)$**. Assume $-2 < \beta_1 < 1$, z is nonzero and $Re(z) \geq 0$. The integral used to define

$$B_0(z, \beta) := (1/z)\int_{\theta>0} d(\theta)\theta^{-(1+\beta)}(1/(1+\theta^2))(1 - e^{-z\theta})$$

converges absolutely. $B_0(z, \beta)$ is continuous in those z and analytic in those $\beta$. Also $B_0(z, \beta)$ is analytic in z for $Re(z) > 0$. Say $-1 < \beta_1 < 1$. $B_0(z, \beta)$ is then continuous at $z = 0$. $B_0(0, \beta) := (\pi/2)/\cos((\pi/2)\beta)$.

**Lemma 3.7 Mellin transform representation of E(u, β, 1).** *Assume* $-2 < \beta_1 <$

*1 and max{0, –(1 + $\beta_1$)} < Re(u) < min{1, 1 - $\beta_1$}.*

$$E(u, \beta, 1) = \int_{v>0} d(v) v^{u-1} B_0(v, \beta),$$

*with absolute convergence for the integral.*

**The Mellin transform representation of F(u, β, 1).**

**Definition and properties of M(z, β).** Set $\beta_1 := \mathrm{Re}(\beta)$. Suppose $-2 < \beta_1 < 1$. Say z is nonzero and Re(z) ≥ 0. Take

$$M(z, \beta) := (2/\pi)\sin(\pi\beta)B_0(z, \beta) + (2/z)(\cos((\pi/2)\beta) - \cos(z - (\pi/2)\beta)).$$

M(z, β) is continuous in z and analytic in β. Also M(z, β) is analytic in z for Re(z) > 0. Say $-1 < \beta_1 < 1$. M(z, β) has a continuous extension to z = 0, with M(0, β) = 0.

**Lemma 3.8 Mellin transform representation of F(u, β, 1).** *Assume -2 < $\beta_1$ < 1 and max{0, –(1 + $\beta_1$)} < Re(u) < min{1, 1 - $\beta_1$}.*

$$F(u, \beta, 1) = \int_{v>0} d(v) v^{u-1} M(v, \beta),$$

*with absolute convergence for the integral.*

Assume $-2 < \beta_1 < 0$. W(0, 1 + β) = -(π/2)/sin((π/2)β). Say Re(z) ≥ 0. $zB_0(z, \beta)$ = -(π/2)/sin((π/2)β) - W(z, 1 + β). So M(z, β) = (-2/z)((1/π)sin(πβ)W(z, 1 + β) + cos(z – (π/2)β)). Suppose Re(z) > 0. Apply the translation relation for W(z, 1 + β). Then

$$M(z, \beta) = (2/z)((1/\Gamma(1 + \beta))(z^{\beta}) + (1/\pi)\sin(\pi\beta)W(z, \beta - 1) - \cos(z - (\pi/2)\beta)).$$

That formula for M(z, β) extends to the case with Re(z) > 0 and $-2 < \beta_1 < 1$ as well as to that with z ≠ 0, Re(z) ≥ 0 and $0 < \beta_1 < 1$.

W(z, 1 + β) will be determined, for certain ranges of z and β, from I(p, z, u) and so from I(p, z), as presented next.

**Definition and properties of I(p, z, u).**

**Claim 3.2** *Assume either (1) or (2) as stated next.*
(1) *Re(p) > 0 and Re(z) > 0. u ≥ 0 or u is not real.*
(2) *0 < Re(p) < 1 and Re(z) ≥ 0. u > 0 or u is not real.*
*Then* $\int_{j>0} d(j)\, )\cdot|j^{p-1}/(1 + uj))e^{-zj}|$ *converges. Set*

$$I(p, z, u) := \int_{j>0} d(j)(j^{p-1}/(1 + uj))e^{-zj}.$$

*I(p, z, u) is analytic in each of p, z, u with Re(z) > 0 and u nonzero. If (1), then I(p, z, u) is continuous at u = 0. When (2) holds, I(p, z, u) is continuous in z.*

**Proof** Say $\pi/2 \le B < \pi$. Assume (*): u is nonzero and $\arg(u) \le B$.
Consider $\int_{0<j<1} d(j)j^{-1+Re(p)}(1/|(1 + uj)|)\, e^{-Re(z)j}$. $e^{-Re(z)j} \le 1$. If $r \ge 0$, then $1/|(1 + ur)| \le 1/\sin(B)$.
Consider $\int_{j\ge 1} d(j)j^{-2+Re(p)}(1/|(j^{-1} + u)|)e^{-Re(z)j}$. Say d > 0. Assume (*) and $|u| \ge d$. If $r \ge 0$, then $1/|(r + u)| \le 1/(d\cdot\sin(B))$.

**The determination of W(z, 1 + β) from I(p, z, u).** Assume either $-1 < \beta_1 < 0$ and $Re(z) \ge 0$, or $\beta_1 < 0$ and Re(z) > 0. Then

$$W(z, 1 + \beta) = \tfrac{1}{2}\textstyle\sum_{\sigma=\pm 1} I(-\beta, z, \sigma i).$$

**The determination of I(p, z, u) from I(p, z/u).**

**Definition of I(p, z).** Set I(p, z) := I(p, z, 1).

Let u > 0. Say either $0 < p_1 < 1$ and $Re(z) \ge 0$, or $p_1 > 0$ and Re(z) > 0. One then has $I(p, z, u) = u^{-p}\cdot I(p, z/u)$. We will extend that equality to certain complex u.

Say u is nonzero and is not a negative real. We use the principal branch for which $|\arg(u)| < \pi$ for $u^\beta$.

**Claim 3.3** *Assume $0 < p_1 < 1$, u is not zero, Re(z) > 0 and $Re(z/u) \ge 0$.
Then (1): $I(p, z, u) = u^{-p}\cdot I(p, z/u)$.
When z is positive, (1) holds for all u with $Re(u) \ge 0$.*

**Proof** If u > 0, then Re(z/u) > 0. (1) holds for u > 0. Each side of the asserted equality (1) is analytic in the nonzero u with Re(z/u) > 0. Hence the equality holds also for those u. Each side of (1) is continuous in the nonzero u with $Re(z/u) \ge 0$. So (1) holds also for nonzero u with Re(z/u) = 0.

**The determination of W(z, 1 + β) from I(-β, ±iz).** Say $-1 < \beta_1 < 0$ and $z \ge 0$. Then $I(-\beta, z, \sigma i) = (\sigma i)^\beta I(-\beta, -\sigma iz)$. So

$$W(z, 1 + \beta) := \tfrac{1}{2}\textstyle\sum_{\sigma=\pm 1} (\sigma i)^\beta I(-\beta, -\sigma iz).$$

I(β, z) is determined in Lemma 4.1.

## §4 The incomplete gamma functions and confluent hypergeometric functions

We now derive and adapt to our needs some classical results related to the incomplete gamma functions and confluent hypergeometric functions. See H. Buchholz [4], Erdélyi et al. [9], L. J. Slater [12-13] and F. G. Tricomi [14].

**Relation between I(p, z) and Γ(1 – p, z).**

**Lemma 4.1** *Assume either (1) or (2) as stated next.*
(1) *$Re(p) > 0$ and $Re(z) > 0$.*
(2) *$0 < Re(p) < 1$ and $Re(z) \geq 0$.*
*Then*

$$(1/\Gamma(p))e^{-z}I(p, z) = \Gamma(1 - p, z).$$

**Proof** Assume (1). $e^{-z}I(p, z) = \int_{j>0} d(j)(j^{p-1}/(1+j))e^{-z(1+j)}$. Then $(d/dz)(1/\Gamma(p))e^{-z}I(p, z) = -e^{-z}z^{-p}$.
Now consider all z with $z > 0$ or z nonreal. Use $|\arg(z)| < \pi$. $(d/dz)\Gamma(q, z) = -e^{-z}z^{q-1}$. If $Re(q) > 0$, then $\Gamma(q, z)$ has a continuous extension to $z = 0$: $\Gamma(q, 0) := \Gamma(q)$. Also $\Gamma(q, z)$ is analytic in q.
Assume (1). Set $\Delta(p, z) := (1/\Gamma(p))e^{-z}I(p, z) - \Gamma(1 - p, z)$. On the half-plane of z with $Re(z) > 0$: $(d/dz)\Delta(p, z) = 0$, hence $\Delta(p, z) = K(p)$. Assume $0 < Re(p) < 1$. $I(p, z)$ and $\Gamma(1 - p, z)$ are continuous for $Re(z) \geq 0$. $I(p, 0) := \pi/\sin(\pi p)$. Therefore $\Delta(p, 0) = 0$. So $\Delta(p, z) = 0$ for $Re(z) \geq 0$.
Assume $Re(z) > 0$. $(1/\Gamma(p))e^{-z}I(p, z)$ is analytic in all p with $Re(p) > 0$. So is $\Delta(p, z)$. Hence $\Delta(p, z) = 0$ for $Re(p) > 0$.

The Mellin transform representation of $F(z, \beta)$ is obtained, in Lemma 5.1 and Theorem 5.2, from the function $H(z, \beta)$ arising from $\varphi(1 + \beta, u)$, as developed next.

**Definition of φ(B, z).** Let $|z|$ and $|Im(B)|$ be bounded. Say $Re(B) \geq r$. $\sum_{k\geq 0} |z^k/\Gamma(B + k)|$ converges uniformly in the aforementioned z, B. This is seen on applying:
(1) Say $\delta > 0$, y is real and $|y| \leq \delta$. If $x \geq 1$, then $1/|\Gamma(x + iy)| \leq (\sinh(\pi\delta)/(\pi\delta))^{1/2}(1/(\Gamma(x)))$.
(2) Stirling's formula for $\Gamma(x)$.
So

$$\varphi(B, z) := \sum_{k\geq 0} z^k/\Gamma(B + k)$$

is entire in each of B, z.

**Definition of $\gamma(\beta, z, *)$.** Assume $\mathrm{Re}(\beta) > 0$. Say $z \geq 0$ or z is nonreal. $z^{-\beta}\cdot\gamma(\beta, z) = \int_{0 < j < 1} d(j) j^{\beta-1} e^{-zj}$. So

$$\gamma(\beta, z, *) := \sum_{k \geq 0} (-z)^k/(k!(\beta + k)\Gamma(\beta))$$

extends $z^{-\beta}\cdot\gamma(\beta, z)/\Gamma(\beta)$ to a function entire in z and β.

**The confluent hypergeometric functions M(a, B, z) and U(a, a, z).**

When B is not zero or a negative integer, set $M(1, B, z) := \Gamma(B)\varphi(B, z)$. Then $M(1, B, z) := \sum_{k \geq 0} z^k/(B)_k$. M(1, B, z) is a special case of the confluent hypergeometric function of the first kind, M(a, B, z).

**Definition of M(a, B, z).**

M(a, B, z) satisfies the equation specified next.

**Definition of K(a, B).** $K(a, B) := z(D^2) + (B - z)D - a$, where $D := d/dz$.

**Kummer's equation**: $K(a, B)(g) = 0$. Initial conditions: $M(a, B, 0) := 1$, $M_z(a, B, 0) := a/B$.

$K(a, B) = D\cdot(z(D - 1) + B - 1) + 1 - a = (zD + a)(D - 1) + (B - a)D$.

**Relation between $\varphi(1 + \beta, z)$ and $\gamma(\beta, z, *)$.**

**Lemma 4.2** *$\varphi(1 + \beta, z) = e^z\cdot\gamma(\beta, z, *)$.*

**Proof** Let $\mathrm{Re}(\beta) > 0$. Set $B = 1 + \beta$. Say z is nonzero and $|\arg(z)| < \pi$. Then

**Eq. (*)**

$$(d/dz)(z^{\beta}\varphi(B, z)) = z^{\beta-1}/\Gamma(\beta) + z^{\beta}\varphi(B, z).$$

So $(d/dz)(e^{-z}(z^{\beta}\varphi(B, z))) = e^{-z}z^{\beta-1}/\Gamma(\beta)$.
$\mathrm{Re}(\beta) > 0$ implies $0^{\beta} = 0$. Thus $\varphi(1 + \beta, z) = e^z\cdot\gamma(\beta, z, *)$.
$\varphi(1 + \beta, z) = e^z\cdot\gamma(\beta, z, *)$ extends to all β, z, since $\varphi(1 + \beta, z)$ and $e^z\cdot\gamma(\beta, z, *)$ are entire in β, z.

The equation (*) amounts to $(zD + \beta - z)\varphi(1 + \beta, z) = 1/\Gamma(\beta)$. $K(1, B) = D\cdot(z(D - 1) + B - 1)$. Therefore $\varphi(1 + \beta, z)$ satisfies Kummer's equation with $a = 1$ and $B = 1 + \beta$.

**Definition of U(a, a, z).** We have seen that $(D - 1)(e^z\Gamma(1 - a, z)) = z^{-a}$. Now $K(a, a) = (zD + a)(D - 1)$. So $K(a, a)(e^z\Gamma(1 - a, z)) = 0$. Kummer's function of the second kind $U(a, a, z)$ is one of the solutions of the latter equation with $U(a, a, 0) = \Gamma(1 - a)$. It turns out that $U(a, a, z) = e^z\Gamma(1 - a, z)$.

**The Laplace representation of φ(1 + β, z).** If $\mathrm{Re}(\beta) > 0$, then

$$\varphi(1 + \beta, z) = (1/\Gamma(\beta))\int_{0 < j < 1} d(j)(1 - j)^{\beta - 1}e^{zj}.$$

**Relation between φ(1 + β, z) and φ(β, z).** Let z be nonzero and $|\arg(z)| < \pi$. Assume $\mathrm{Re}(\beta) > 0$. Then $\gamma(1 + \beta, z) = -e^{-z}z^\beta + \beta\gamma(\beta, z)$. Therefore, for all β, z, $z\varphi(1 + \beta, z) = -1/\Gamma(\beta) + \varphi(\beta, z)$.

**Relation between I(-β, z) and φ(1 + β, z).** Let z be nonzero and $|\arg(z)| < \pi$. $\Gamma(1 + \beta, z) = e^{-z}z^\beta + \beta\Gamma(\beta, z)$. Assume $\mathrm{Re}(\beta) > 0$. $\Gamma(\beta, z) = \Gamma(\beta) - \gamma(\beta, z)$. So $z^{-\beta}\Gamma(\beta, z)/\Gamma(\beta) = z^{-\beta} - \gamma(\beta, z, *)$. The latter relation extends to all β. Thus $z^{-\beta}e^z\Gamma(1 + \beta, z)/\Gamma(1 + \beta) = 1/\Gamma(1 + \beta) - \varphi(1 + \beta, z) + z^{-\beta}e^z$.
Say $-1 < \mathrm{Re}(\beta) < 0$, $\mathrm{Re}(z) \geq 0$. $z^{-\beta}I(-\beta, z) = \Gamma(-\beta)z^{-\beta}e^z\Gamma(1 + \beta, z)$. So $z^{-\beta}I(-\beta, z) = -(\pi/\sin(\pi\beta))(1/\Gamma(1 + \beta) - \varphi(1 + \beta, z) + z^{-\beta}e^z)$.

## §5 H(z, β) and the Mellin transform representation of F(z, β).

**Definition of H(z, β).** Set

$$H(z, \beta) := 2/\Gamma(1 + \beta) - \sum_{\sigma = \pm 1} \varphi(1 + \beta, \sigma iz)).$$

Let $\beta_1 := \mathrm{Re}(\beta)$. Assume $-1 < \beta_1 < 0$. Take $z > 0$.

$$W(z, 1 + \beta) := z^\beta \cdot \tfrac{1}{2}\sum_{\sigma = \pm 1} (-\sigma iz)^{-\beta}I(-\beta, -\sigma iz).$$

**The determination of W(z, 1 + β) from H(z, β).**

$$W(z, 1 + \beta) = -(\pi/\sin(\pi\beta))(z^\beta \cdot \tfrac{1}{2}H(z, \beta) + \cos(z - (\pi/2)\beta)).$$

The latter expression for $W(z, 1 + \beta)$ extends to the half-plane of z with $\mathrm{Re}(z) \geq 0$.

Note that $H(z, \beta)$ is entire in each of z, β. $H(-z, \beta) = H(z, \beta)$. $(H(z, \beta))^* = H(z^*, \beta^*)$.

**Definition of Ω(r, δ).** Assume r is real and $\delta \geq 0$. Let $\Omega(r, \delta)$ be the set of all $\beta = \beta_1 + i\beta_2$ with $\beta_1 \geq r$ and $-\delta \leq \beta_2 \leq \delta$.

Say $\varepsilon \geq 0$. There exists a $K_0(r, \delta, \varepsilon) > 0$ such that $|H(z, \beta)| \leq K_0(r, \delta, \varepsilon)(|z|^2)$, when $\beta$ is in $\Omega(r, \delta)$ and z is in the disc $|z| \leq \varepsilon$.

**The determination of M(z, β) from W(z, 1 + β).** Say $-2 < \beta_1 < 0$. Let $Re(z) \geq 0$.

$$M(z, \beta) = (-2/z)((1/\pi)\sin(\pi\beta)W(z, 1 + \beta) + \cos(z - (\pi/2)\beta)).$$

**The determination of M(z, β) from H(z, β).** Assume $-1 < \beta_1 < 0$. Then

$$M(z, \beta) = z^{-1+\beta} \cdot H(z, \beta).$$

The latter equality extends to the case $-2 < \beta_1 < 1$, $z \neq 0$ and $|\arg(z)| \leq \pi/2$. If $-1 < \beta_1 < 1$ that relation extends by continuity to $z = 0$.

Assume $-2 < \beta_1 < 1$ and $\max\{0, -(1 + \beta_1)\} < Re(u) < \min\{1, 1 - \beta_1\}$.

$$F(u, \beta, 1) = \int_{v>0} d(v)v^{u-1}v^{-1+\beta} \cdot H(v, \beta).$$

Set $z = 1 - (u + \beta)$. $F(z, \beta) = F(u, \beta, 1)$.

**Lemma 5.1** *Assume $-2 < \beta_1 < 1$ and $\max\{0, -\beta_1\} < Re(z) < \min\{2, 1 - \beta_1\}$. Then* **Mellin transform representation of F(z, β)**

$$F(z, \beta) = \int_{v>0} d(v)v^{z-1}H(1/v, \beta).$$

The previous Lemma 5.1 will be improved upon in the Theorem 5.2 on the Mellin transform representation of F(z, β), after obtaining properties of H(z, β) to be used for that purpose.

**Integrability of |F(z, β)| on vertical lines.** $|F(z, \beta)|$ is square integrable on $x + i\mathbf{R}$, for any $x > 0$, provided $\beta_1 \geq 0$.

$H(z, 0) := 2(1 - \cos(z))$. Given $L \geq 0$, restrict z to the horizontal strip of all z′ with $|Im(z')| \leq L$. Say $z = z_1 + iz_2$ with $z_1$, $z_2$ real. $|H(z, 0)| \leq 2(1 + \cosh(L))$, since $\frac{1}{2}(1 - \cos(z)) = (\sin(z/2))^2$ and $|\sin(z/2)|^2 = \frac{1}{2}(\cosh(z_2) - \cos(z_1))$.

**Claim 5.1 Representation of H(z, β).** *If $\beta_1 > 0$, then*

$$\tfrac{1}{2}H(z, \beta) = (1/\Gamma(\beta))\int_{0<j<1} d(j)(1-j)^{\beta-1}(1 - \cos(zj)).$$

**Proof** If $\beta_1 > 0$, then $\varphi(1 + \beta, z) = (1/\Gamma(\beta))\int_{0<j<1} d(j)(1-j)^{\beta-1}e^z$.

Set m := min {$\Gamma(x)$: $1 < x < 2$}.

**Corollary 5.1 Uniform boundedness.** *Assume $\varepsilon > 0$, $\delta > 0$ and $L \geq 0$. $|H(z, \beta)|$ is bounded for $\beta$ in $\Omega(\varepsilon, \delta)$ and z in the strip $|z_2| \leq L$: $|H(z, \beta)| \leq K(\varepsilon, \delta, L)$.*

**Proof** If $\varepsilon = 1$, then $K(1, \delta, L) := 2(1 + \cosh(L))(\sinh(\pi\delta)/(\pi\delta))^{½}(1/m)$ works. If $0 < \varepsilon < 1$, take $K(\varepsilon, \delta, L) := (1 + \delta/\varepsilon)K(1, \delta, L)$ .

**Corollary 5.2 Positivity.** *Assume $\beta \geq 0$. Say z is real and not zero. Then $H(z, \beta) > 0$.*

**Translation relation for H(z, β).**

**Definition of A(z, β, n).** Say n is an integer and $z \neq 0$. Set $A(z, \beta, n) := \sum_{0 \leq k \leq n} (1/\Gamma(\beta + 1 - 2k))(-1/(z^2))^k$.

**Claim 5.2** *Assume z is not zero.*
(1) *$½H(z, \beta) = 1/\Gamma(\beta + 1) - (1/z^2)½H(z, \beta - 2)$.*
(2) *Say n is an integer $\geq 0$. $½H(z, \beta) := A(z, \beta, n - 1) + (-1/(z^2))^n \cdot ½H(z, \beta - 2n)$.*

**Proof** $½H(z, \beta) := 1/\Gamma(1 + \beta) - ½\sum_{\sigma = \pm1} \varphi(1 + \beta, \sigma iz)$. Now $z\varphi(\beta + 1, z) = -1/\Gamma(\beta) + \varphi(\beta, z)$ for all β, z. Thus $z\sum_{\sigma = \pm1} \varphi(\beta + 1, \sigma z) = \sum_{\sigma = \pm1} \sigma\varphi(1 + \beta, \sigma z)$. Therefore $z^2 \sum_{\sigma = \pm1} \varphi(\beta + 2, \sigma z) = - 2/\Gamma(\beta) + \sum_{\sigma = \pm1} \varphi(\beta, \sigma z)$. Thus (1) holds. Hence so does (2).

**The determination of H(z, β) from W(z, β - 1).**

**Claim 5.3** *Say $z \neq 0$, $|arg(z)| < \pi/2$ and $\beta_1 < 2$. Then*
(1) *$½H(z, \beta) = 1/\Gamma(\beta + 1) - (z^{-\beta})\cos(z - (\pi/2)\beta) + z^{-2}(\sin(\pi\beta)/\pi)(z^{2-\beta}W(z, \beta - 1))$.*

**Proof** If $-1 < \beta_1 < 0$ and $Re(z) \geq 0$, then $H(z, \beta) = 2(z^{-\beta})(-(\sin(\pi\beta)/\pi)W(z, 1 + \beta) - \cos(z - (\pi/2)\beta))$. If $\beta_1 < 0$ and $Re(z) > 0$, then $W(z, 1 + \beta) = \Gamma(-\beta)z^\beta - W(z, \beta - 1)$. Thus, for $Re(z) > 0$ and $-1 < \beta_1 < 0$, $½H(z, \beta) = (z^{-\beta})((\sin(\pi(-\beta))/\pi)(\Gamma(-\beta)z^\beta - W(z, \beta - 1)) - \cos(z - (\pi/2)\beta))$. So (1) holds. That equality extends to $Re(z) > 0$ and $\beta_1 < 2$ since each side is analytic in β there.

**The representation of W(z, β) via R(z, β).**

If $\beta_1 < 1$ and $z > 0$, then $z^{1-\beta}W(z, \beta) = \int_{j>0} d(j)j^{-\beta}(1/(1 + (j/z)^2))e^{-j}$. We generalize that equalty to certain complex z next.

**Claim 5.4** *Say $K \leq \beta_1 < 1$. Let B with $0 \leq B < \pi/2$ be given. Restrict z to the closed sector: Either $z \neq 0$ and $|arg(z)| \leq B$, or $z = 0$.*

$$\int_{j>0} d(j)\ |j^{-\beta}(1/(1+(zj)^2))e^{-j}| \le (\cos(B))^{-2}\Gamma(1-K).$$

**Proof** $|1 \pm (iz)j| \ge \cos(B)$. So $|j^{-\beta}(1/(1+(zj)^2))| \le j^{-K}(\cos(B))^{-2}$.

**Definition of R(z, β).** Let $\mathrm{Re}(z) \neq 0$ and $\beta_1 < 1$.
Set $R(z, \beta) := \int_{j>0} d(j) j^{-\beta}(1/(1+(j/z)^2))e^{-j}$.

$R(-z, \beta) = R(z, \beta)$. Say $z = |z|e^{i\theta}$. $|R(z, \beta)| \le (\cos(\theta))^{-2}\Gamma(1-\beta_1)$.

**Claim 5.5** *Say $\mathrm{Re}(z) > 0$ and $\beta_1 < 1$. $R(z, \beta)$ is analytic in each of z, β.*

**Claim 5.6** *Assume $\beta_1 < 1$.*
(1) *$z^{1-\beta}W(z, \beta) = R(z, \beta)$, for all nonzero z with $-\pi/2 < \arg(z) < \pi/2$.*
(2) *Given $L \ge 0$, restrict z to the strip $|z| \le L$.*
   (i) *Say $\mathrm{Re}(z) \ge x > 0$. $|R(z, \beta)| \le (1 + (L^2/x^2))\Gamma(1-\beta_1)$.*
   (ii) *$\lim_{\mathrm{Re}(z)\to\infty} R(z, \beta) = \Gamma(1-\beta)$.*

**Proof of (1).** $z^{1-\beta}W(z, \beta) = R(z, \beta)$ extends from the ray $z > 0$ to the half-plane of z with $\mathrm{Re}(z) > 0$, since $z^{1-\beta}W(z, \beta)$ and $R(z, \beta)$ are each analytic there.

**The boundedness of H(z, β).**

**Lemma 5.2** *Assume $z \neq 0$ and $|\arg(z)| < \pi/2$.*
(1) *Say $\beta_1 < 2$.*

$$\tfrac{1}{2}H(z, \beta) = 1/\Gamma(\beta+1) - (z^{-\beta})\cos(z - (\pi/2)\beta) + z^{-2}(\sin(\pi\beta)/\pi)R(z, \beta-1).$$

(2) *Let n be an integer $\ge 0$. Assume $\beta_1 < 2(n+1)$.*

$$\tfrac{1}{2}H(z, \beta) =$$

$$A(z, \beta, n) - (z^{-\beta})\cos(z - (\pi/2)\beta) - (-1/(z^2))^{n+1}(\sin(\pi\beta)/\pi)R(z, \beta-(2n+1)).$$

*Here*

$$|R(z, \beta-(2n+1))| \le (\cos(\theta))^{-2}\Gamma(2(n+1)-\beta_1),$$

*with $\theta := \arg(z)$.*

**Corollary 5.3** *Say $1 \le \varepsilon < 2$, $\delta > 0$. Let β vary over the rectangle $0 \le \beta_1 \le \varepsilon$, $|\beta_2| \le \delta$. Suppose $L \ge 0$. Let z vary over the strip $|\mathrm{Im}(z)| \le L$. Then $|H(z, \beta)|$ is bounded.*

*On the ray $v \geq 1$, $|H(v, \beta)| \leq K_1(\varepsilon, \delta)$, with $K_1(\varepsilon, \delta) := 2((\sinh(\pi\delta)/(\pi\delta))^{½}(1/m) + 1 + (½(1 + \cosh(2\pi\delta)))^{½}\Gamma(\varepsilon)/\pi)$.*

**Proof** $H(z, \beta)$ is continuous in all $z$, $\beta$. So we may assume $|z| \geq r > 0$. The corollary follows from the previous Lemma 5.2 (1) and Claim 5.6 (2) (i).

**Theorem 5.1** *Given $L \geq 0$, restrict $z$ to the strip $|Im(z)| \leq L$.*
*If $\beta_1 > 0$, $\beta' \to \beta$ and $|Re(z)|$ tends to $\infty$, then $½H(z, \beta')$ converges to $1/\Gamma(1 + \beta)$.*
*If $\beta$ varies over $\Omega(0, \delta)$ and $z$ over the strip, then $|H(z, \beta)|$ is bounded.*

**Theorem 5.2 Mellin transform representation of F(z, β).**
*Suppose $\beta_1 \geq 0$ and $0 < Re(z) < 2$.*
(1) $\int_{v > 0} d(v)|v^{z-1}H(1/v, \beta)|$ *is finite.*
(2) $F(z, \beta) = \int_{v > 0} d(v)v^{z-1}H(1/v, \beta)$.

**Proof of (1).**

Claim 1 *Say $Re(z) < 2$ and $\beta$ is arbitrary. $\int_{v \geq 1} d(v)v^{Re(z)-1}|H(1/v, \beta)|$ is finite. Also $\int_{v \geq 1} d(v)v^{z-1}H(1/v, \beta)$ is analytic in those $z$ and $\beta$.*

Proof of Claim 1. $|H(u, \beta)| \leq K_0(r, \delta, 1)(|u|^2)$, when $\beta$ varies over $\Omega(r, \delta)$ and $|u| \leq 1$. $v^{Re(z)} \leq v^x$, when $v \geq 1$ and $Re(z) \leq x$. Take $x < 2$. $\int_{v \geq 1} d(v)v^{x-2}$ is finite. So the integral $\int_{v \geq 1} d(v)v^{Re(z)-1}|H(1/v, \beta)|$ converges uniformly on any compact set of pairs $z$, $\beta$ with $Re(z) < 2$.

Claim 2 *Say $Re(z) > 0$ and $\beta_1 \geq 0$. $\int_{0 < v < 1} d(v)v^{Re(z)-1}|H(1/v, \beta)|$ is finite. Also for those $z$ and $\beta$, $\int_{0 < v < 1} d(v)v^{z-1}H(1/v, \beta)$ is analytic in $z$, continuous in $\beta$ when $\beta_1 \geq 0$ and analytic in $\beta$ when $\beta_1 > 0$.*

Proof of Claim 2. Let $\beta$ vary over $\Omega(0, \delta)$ and $v$ over $(0, 1]$. Then $|H(1/v, \beta)|$ is bounded. Also $Re(z) \geq x > 0$ implies $v^{Re(z)} \leq v^x$. The integral $\int_{0 < v < 1} d(v)v^{Re(z)-1}|H(1/v, \beta)|$ converges uniformly on any compact set of pairs $z$, $\beta$ with $Re(z) > 0$ and $Re(\beta) \geq 0$.

**Proof of (2).** If $-2 < \beta_1 < 1$ and $\max\{0, -\beta_1\} < Re(z) < \min\{2, 1 - \beta_1\}$, then (2) is valid by Lemma 5.1. Now assume $\beta_1 \geq 0$. Each side of (2) is analytic in the $z$ with $0 < Re(z) < 2$. Fix $\beta$ with $0 \leq \beta_1 < 1$. Then (2) holds if $0 < Re(z) < 1 - \beta_1$. Therefore (2) holds for all $z$ with $0 < Re(z) < 2$. Fix $z$ with $0 < Re(z) < 2$. Then (2) holds for all $\beta$ with $0 \leq \beta_1 < 1$. Each side of (2) is analytic in the $\beta$ with $\beta_1 > 0$. So (2) holds for all $\beta$ with $\beta_1 \geq 0$.

Suppose $\beta_1 > -3$, $m$ is a positive integer and $z$ is complex. $\int_{0 < j < 1} d(j)|j^{\beta}(1-j)^{m-1}H(jz, \beta)|$ is finite, since $H(u, \beta) = O(u^2)$ for $u$ near to 0.

**Definition of $G(z, \beta, m)$ and $H(z, \beta, 2w)$.** Set $G(z, \beta, m) := (1/((m - 1)!))\int_{0 < j < 1} d(j) j^{\beta}(1- j)^{m-1} H(jz, \beta)$. Let $G(z, \beta, 0) := H(z, \beta)$, for all z, β. Say w is a nonnegative integer. Set $H(z, \beta, 2w) := z^{2w} G(z, \beta, 2w)$. $H(z, \beta, 0) = H(z, \beta)$.

**Lemma 5.3**
(1) *Assume β is in Ω(r, δ).*
(i) *If $r > -1$ and z is in the strip $|Im(z)| \leq L$, then $|G(z, \beta, m)| \leq K(r, \delta, L)/(1 + r)_m$.*
(ii) *If $r > -3$ and z is in the disc $|z| \leq \varepsilon$, then $|G(z, \beta, m)| \leq (K_0(r, \delta, \varepsilon)/(3 + r)_m)(|z|^2)$.*
(2) *Say $\beta_1 > -3$.*
(i) *Fix β, m, z. Set $t(j, k) := (1/((m - 1)!)) j^{\beta + 2k}(1- j)^{m-1}(-1)^{k-1}(z^{2k})/\Gamma(1 + \beta + 2k)$. $\sum_{k \geq 1} \int_{0 < j < 1} d(j)|t(j, k)|$ is finite.*
(ii) *$\frac{1}{2}G(z, \beta, m) = (-1)\sum_{k \geq 1} (-1)^k(z^{2k})/\Gamma(1 + \beta + m + 2k)$. Thus $G(z, \beta, m)$ extends to an entire function of each of z, β + m.*
(3) *Say w is a nonnegative integer.*

$$\tfrac{1}{2}((-1)^w H(z, \beta, 2w)) = (-1)\sum_{k \geq w + 1} (-1)^k(z^{2k})/\Gamma(1 + \beta + 2k)$$

$$= \tfrac{1}{2}H(z, \beta) + \sum_{1 \leq k \leq w} (-1)^k(z^{2k})/\Gamma(1 + \beta + 2k).$$

(4) *There is a $K(\beta, w) > 0$ such that $|H(r, \beta, 2w)| \leq K(\beta, w)g(2w, 2(w + 1)(r)$, for all $r > 0$.*

**Proof of (2) (i).** $\frac{1}{2}H(z, \beta) = (-1)\sum_{k \geq 1} (-1)^k(z^{2k})/\Gamma(1 + \beta + 2k)$. Thus $\frac{1}{2}G(z, \beta, m) = \int_{0 < j < 1} d(j)\sum_{k \geq 1} t(j, k)$. Say $r_1 > -3$ and β is in Ω(r, δ). $\sum_{k \geq 1} \int_{0 < j < 1} d(j)|t(j, k)| \leq (1 + \delta/(r + 3))(\sinh(\pi\delta)/(\pi\delta))^{1/2}(-\frac{1}{2}H(i|z|, \beta_1 + m))$.
**Proof of (2) (ii).** $\int_{0 < j < 1} d(j)\sum_{k \geq 1} t(j, k) = \sum_{k \geq 1} \int_{0 < j < 1} d(j)t(j, k)$.

**Theorem 5.3 Mellin transform representation of $F(z, \beta)$.**
*Suppose $\beta_1 \geq 0$, $w \geq 0$ and $2w < Re(z) < 2(w + 1)$.*
(1) *$\int_{v > 0} d(v)|v^{z-1}H(1/v, \beta, 2w)|$ is finite.*
(2) *$(-1)^w F(z, \beta) = \int_{v > 0} d(v)v^{z-1}H(1/v, \beta, 2w)$.*

**Proof** Say $w \geq 1$. Set $m = 2w$. Let $Z = z - m$. Apply the translation relation for $F(z, \beta)$ and then the Mellin transform representation of $F(Z, \beta)$ established in Theorem 5.2.

$$(-1)^w F(z, \beta) = ((1/((m - 1)!))\int_{0 < v \leq 1} d(v)\, v^{Z-1}v^{1+\beta}(1- v)^{m-1})\cdot\int_{v > 0} d(v)v^{Z-1}H(1/v, \beta).$$

**Lemma 5.4 Positivity.** *Say $w \geq 0$.*

(1) *Assume β > -2. H(r, β, 2w) strictly increases with r on [0, 1].*
(2) *Assume β ≥ 0. H(r, β, 2w) > 0, for r a nonzero real. H(0, β, 2w) = 0.*

**Proof of (1).** $\frac{1}{2}H(z, \beta, 2w) = \sum_k (a_k \cdot u^k - a_{k+1} \cdot u^{k+1})$, with $k \geq w + 1$, $k \equiv (w + 1) \bmod(2)$, $u = r^2$ and $a_k = 1/\Gamma(1 + \beta + 2k)$. $(d/du)(a_k \cdot u^k - a_{k+1} \cdot u^{k+1}) \geq k \cdot a_{k+1} u^{k-1}(a_k/a_{k+1} - (1 + 1/k))$. Now $a_k/a_{k+1} = (1 + \beta + 2k)(1 + \beta + 2k + 1)$. So $a_k/a_{k+1} > 2 \geq 1 + 1/k$.

**Proof of (2).** Say $\beta \geq 0$. $H(r, \beta) > 0$ for r a nonzero real, by the previous corollary on positivity. Let $m \geq 1$. $G(r, \beta, m) := (1/((m - 1)!))\int_{0<j<1} d(j) j^\beta (1 - j)^{m-1} H(jr, \beta) > 0$. Say $w \geq 1$. $H(r, \beta, 2w) := r^{2w} G(r, \beta, 2w) > 0$.

## §6 The Mellin transform representation of f(s, β) := 1/(sin(πs/4)·2ξ(2β + s)).

**Definitions of $n_0(s, \beta)$, $f_0(s, \beta)$, $b(s, \beta)$, $T_0(z, \beta)$ and $T_0(z, \beta, 4w)$.**

$$n_0(s, \beta) := \sin(\pi s/4) l(2\beta + s).\ f_0(s, \beta) := 1/n_0(s, \beta).\ b(s, \beta) := n_0(s, \beta)(2\beta + s - 1).$$

$$T_0(z, \beta) := (-1/(\pi^{1-\beta}\Gamma(1 + \beta)))\sum_{k\geq 1} (\pi z)^{2k}/((1 + \beta)_{2k}(2k + \beta - \tfrac{1}{2})).$$

$(1/\Gamma(\frac{1}{2}(3/2 + \beta)))T_0(z, \beta)$ is an entire function of each of z, β.

$$T_0(z, \beta) =$$

$$(-1/(\pi^{1-\beta}\Gamma(1 + \beta)))(1/(\tfrac{1}{2} - \beta)(1 - (\tfrac{1}{2})\sum_{\sigma=\pm 1} {}_2F_2(1, \beta - \tfrac{1}{2}; 1 + \beta, \beta + \tfrac{1}{2}; \sigma\pi z)).$$

See Eric W. Weisstein [15] for ${}_pF_q$.

Say $w \geq 0$.

$$T_0(z, \beta, 4w) := ((-1)^{w+1}/(\pi^{1-\beta}\Gamma(1 + \beta)))\sum_{k\geq w+1} (\pi z)^{2k}/((1 + \beta)_{2k}(2k + \beta - \tfrac{1}{2})).$$

Thus $T_0(z, \beta, 0) = T_0(z, \beta)$.

**Definitions of n(s, β) and f(s, β).**

$$n(s, \beta) := \sin(\pi s/4) \cdot 2\xi(2\beta + s) = b(s, \beta)\zeta(2\beta + s).\ f(s, \beta) := 1/n(s, \beta).$$

Note that n(s) = n(s, ¼) and f(s) = f(s, ¼).

**Relations used to determine the Mellin transform representation of f(s, β).** The following relations will be used to derive the Mellin transform representation of f(s, β) from that of F(z, β). That involves successively obtaining the Mellin transform representation of $f_0(s, \beta)$, $1/b(s, \beta)$, $f(s, \beta)$ from that of $F(s/2, \beta)$, $f_0(s, \beta)$, $1/b(s, \beta)$ respectively.

$$f_0(s, \beta) = \pi^{-1+\beta}(½)\pi^{s/2}F(s/2, \beta).$$

$$1/b(s, \beta) := (1/(s + 2\beta - 1))f_0(s, \beta).$$

$$f(s, \beta) := (1/\zeta(2\beta + s))(1/b(s, \beta)).$$

**Integrability of $1/|b(s, \beta)|$ on vertical lines.** Say $s = x + it$, with x, t real. $1/|b(s, \beta)| \sim O(|t|^{-u})$, with $u = 1 + (1 + 2\beta_1 + x)/2$ and $|t|$ large. $1/|b(x + it, \beta)|$ is integrable (respectively: square integrable) in t, when $\beta_1 \geq 0$ (respectively: $\beta_1 \geq -3/2$), $x > 0$ and x is not a multiple of 4.

**Corollary 6.1** *Say $\beta_1 \geq 0$ and $w \geq 0$. The integrals of (1), (2) below are absolutely convergent.*

**The Mellin transform representation of $f_0(s, \beta)$.**

(1) *Say s is on $V_{4w}$: $4w < Re(s) < 4(w + 1)$.*

$$(-1)^w f_0(s, \beta) = \int_{v>0} d(v)v^{s-1}H_0(v^{-2}, \beta, 4w), \text{ with } H_0(z, \beta, 4w) := \pi^{-1+\beta}H(\pi z, \beta, 2w).$$

**The Mellin transform representation of $1/b(s, \beta)$.**

(2) (i) *Fix $\beta$. Set $m = \max\{0, 1 - 2\beta_1\}$. If $m < Re(s) < 4$, then*

$$1/b(s, \beta) = \int_{v>0} d(v)v^{s-1}T_0(i\cdot v^{-2}, \beta).$$

*$|T_0(i\cdot r^2, \beta)| \leq Kg(m, 4)(r)$, for some $K \geq 0$.*

(ii) *Let $w \geq 1$ and s be on $V_{4w}$.*

$$(-1)^w/b(s, \beta) = \int_{v>0} d(v)v^{s-1}T_0(i\cdot v^{-2}, \beta, 4w).$$

*There is a $K \geq 0$ such that $|T_0(i\cdot r^2, \beta, 4w)| \leq Kg(4w, 4(w + 1))(r)$.*

**Proof of (1).**

$$f_0(s, \beta) = \pi^{-1+\beta}(½)\pi^{s/2}F(s/2, \beta).$$

Apply the previous Theorem 5.3.

**Proof of (2).**

$$(-1)^w/b(s, \beta) = (1/(s - 2\alpha))((-1)^w f_0(s, \beta)), \text{ with } \alpha = ½ - \beta.$$

Part (1) gives

$(-1)^w f_0(s, \beta) = \int_{v>0} d(v)v^{s-1}h(1/v)$, with $h(r) = \pi^{-1+\beta} H(\pi\cdot r^2, \beta, 2w)$, for positive r.

Apply Lemma 5.3 (4). Next employ Lemma 3.2 (1), (2), with j = 4w and q = 4(w + 1). Note that $h^{<2\alpha>}(J) = h_0^{<\alpha>}(J^2)$, with $h_0(r) = \pi^{-1+\beta}\cdot½H(\pi r, \beta, 2w)$. The Taylor series expansion of H(z, β, 2w) leads to that for $h^{<2\alpha>}(r)$ as in Lemma 3.2 (3). The result is that $h^{<2\alpha>}(r) = T_0(i\cdot r^2, \beta, 4w)$. $|h^{<2\alpha>}(r)| \leq Kg(\max\{Re(2\alpha), 4w\}, 4(w + 1))(r)$, by Lemma 3.2 (1).

**Definition of c(z, β).** c(z, β) := 1/(n′(z, β)), with ′ = d/ds and z such that n′(z, β) is not zero.

When k is a nonnegative integer,

$$c(4k, \beta) = (2/\pi)(-1)^k(1/\xi(2\beta + 4k)) =$$

$$(-1)^k(\pi^{2k})/(\pi^{1-\beta}\Gamma(1 + \beta + 2k)(2k + \beta - ½)\zeta(2\beta + 4k)).$$

The zeros of ξ(s) are precisely the nonreal zeros of ζ(s). We avoid the occurrence of poles in any of the c(4k, β) by restricting β. Assume that either Re(β) ≥ -3/2 or there exists a positive integer k for which -2k – 3/2 ≤ Re(β) ≤ -2k . Let D be the set of such β.

Say w is a nonnegative integer.

**Definition of $P_{4w}(z, \beta)$.**

$$P_{4w}(z, \beta) :=$$

$$((-1)^{w+1}/(\pi^{1-\beta}\Gamma(1 + \beta)))\sum_{k\geq w+1} (-1)^k z^{2k}/((1 + \beta)_{2k}(2k + \beta - ½)\zeta(2\beta + 4k)).$$

$P_{4w}(\pi\cdot z, \beta) = (-1)^{w+1}\sum_{k\geq w+1} c(4k, \beta)\cdot z^{2k}$. Each $P_{4w}(z, \beta)$ is an entire function of z. Also $P_{4w}(z, \beta)$ is analytic in β on the interior of D and continuous in β on D.

$P_{4w}(0, \beta) = 0$. $P_{4w}(z, \beta) = O(z^{2(w+1)})$, for z near to 0.

The Main unconditional theorem of Part I, §3, is the special case with β = ¼ of the next theorem.

**The Mellin transform representation of f(s, β).**

**Unconditional theorem 6.1** *Let $\beta_1 := Re(\beta) \geq 0$. Fix β. Set $m = \max\{0, 1 - 2\beta_1\}$.*

*The integrals and series in (1), (2) below are absolutely convergent.*

**Mellin transform representation of f(s, β).**

(1) *If $m < Re(s) < 4$, then*

$$f(s, \beta) = \int_{v>0} d(v) v^{s-1} P_0(\pi \cdot v^{-2}, \beta).$$

(2) *Say $w \geq 1$. If $4w < Re(s) < 4(w + 1)$, then*

$$(-1)^w f(s, \beta) = \int_{v>0} d(v) v^{s-1} P_{4w}(\pi \cdot v^{-2}, \beta).$$

**The determination of $P_{4w}$ from $T_0$.**

(3) *Let $w \geq 0$. $P_{4w}(\pi \cdot z, \beta) = \sum_{n \geq 1} (\mu(n)/(n^{2\beta}))T_0(iz/(n^2), \beta, 4w)$.*

**The order of $P_0(r, \beta)$ in r.**

(4) *Say $r > 0$ and r is large. $|P_0(r, \beta)| \sim O(r^{m/2})$. Say $w \geq 1$. $|P_{4w}(r, \beta)| \sim O(r^{2w})$.*

**Proof**

$$(-1)^w f(s, \beta) := (1/\zeta(2\beta + s))((-1)^w/b(s, \beta)).$$

Apply (2) of the previous Corollary 6.1, with α and h(r) as in its proof. Then

$$(-1)^w/b(s, \beta) = \int_{v>0} d(v) v^{s-1} T(1/v), \text{ with } T(r) = h^{<2\alpha>}(r) = T_0(i \cdot r^2, \beta, 4w).$$

Also $|T(r)| \leq Kg(\max\{Re(2\alpha), 4w\}, 4(w + 1))(r)$, for some $K \geq 0$. Next use Lemma 2.4 and Corollary 3.1, with $p = 2\beta$.

$$(1/\zeta(2\beta + s))((-1)^w/b(s, \beta)) = \int_{v>0} (dv) v^{s-1} \omega_{T, 2\beta}(1/v).$$

$T(z) = T_0(i \cdot z^2, \beta, 4w)$ provides the analytic extension of T(z) from the positive reals to the complex plane. Now apply Lemma 2.1, with $p = 2\beta$, $m = 4(w + 1)$ and $q = 1$. The result is that $\omega_{T, 2\beta}(z) = P_{4w}(\pi \cdot z, \beta)$.

**Results when β = ¼ .**

Take β = ¼ in the Unconditional theorem 6.1 to obtain Theorem 3.2 (ii) and Main unconditional theorem (1), (4) (i) presented in Part I, §3.

**Theorem 3.2** (ii) *Say $\varepsilon > 0$. Then for $v \geq 0$, $|P_0(v)| \sim O(v^{1/4+\varepsilon})$, as $v \to \infty$.*

**Main unconditional theorem**

**(1)** ***On the strip $V_0'$ of s with ½ < Re(s) < 4: f(s) is analytic and***

$$f(s) := 1/(\sin(\pi s/4)\cdot 2\xi(\tfrac{1}{2} + s)) = \int_R d(y)e^{sy}P_0(\pi e^{-2y}), \text{ with}$$

$$P_0(z) := (4/(\pi^{3/4}\Gamma(\tfrac{1}{4})))\cdot(-1)\sum_{k\geq 1}(-(z^2))^k/((5/4)_{2k}(2k-\tfrac{1}{4})\zeta(\tfrac{1}{2} + 4k)) \text{ entire.}$$

**(4) (i)** ***Say w ≥ 1. On $V_{4w}$, f(s) is analytic and***

$$(-1)^w f(s) = \int_R d(y)e^{sy}P_{4w}(\pi e^{-2y}),$$

***with***

$$P_{4w}(z) := (-1)^w(P_0(z) + \sum_{1\leq k\leq w} \tilde{c}(4k)(-(z^2))^k) = (-1)^{w+1}\sum_{k\geq w+1}\tilde{c}(4k)(-(z^2))^k$$

***entire.***

$\beta = \frac{1}{4}$ in the Unconditional theorem 6.1 (3) gives the generalization to $w \geq 0$ and complex z of the conjecture of §2:

$$P_{4w}(\pi\cdot z) = \sum_{n\geq 1} (\mu(n)/(n^{\frac{1}{2}}))T_0(iz/(n^2), 4w).$$

**Positivity when β and w are nonnegative.**

**Lemma 6.1 Positivity** *Suppose β ≥ 0 and w ≥ 0.*
(1) *If $P_{4w}(r, \beta)$ is positive for all r > 0, then so is $T_0(ir, \beta, 4w)$.*
(2) *$T_0(ir, \beta, 4w)$ is positive for all r > 0.*

**Proof of (1).** Say $\beta_1 := \mathrm{Re}(\beta) \geq 0$. $(-1)^w/b(s, \beta) = \zeta(2\beta + s)((-1)^w f(s, \beta))$. $(-1)^w/b(s, \beta) = \int_{v>0} d(v)v^{s-1}T_0(i\cdot v^{-2}, \beta, 4w)$, as in (2) of the previous corollary. $(-1)^w f(s, \beta) = \int_{v>0} d(v)v^{s-1}P_{4w}(\pi\cdot v^{-2}, \beta)$ as in (1), (2) of the previous theorem. So $\int_{v>0} d(v)v^{s-1}T_0(i\cdot v^{-2}, \beta, 4w) = \zeta(2\beta + s)\int_{v>0} d(v)v^{s-1}P_{4w}(\pi\cdot v^{-2}, \beta)$.

Fix β, w. Let $T(r) = P_{4w}(\pi\cdot r^2, \beta)$. Part (4) of the previous theorem yields $|T(r)| \leq Kg(r, j, 4(w + 1))$, with $j = \max\{0, 1 - 2\beta_1\}$, when w = 0, and j = 4w when w ≥ 1. Apply Lemma 2.4. $\zeta(2\beta + s)\int_{v>0} d(v)v^{s-1}T(1/v) = \int_{v>0} (dv)v^{s-1}\theta_{T,2\beta}(1/v)$, with $\theta_{T,2\beta}(r) := \sum_{n\geq 1} (n^{-2\beta})P_{4w}(\pi\cdot(r/n)^2, \beta)$.
So $\int_{v>0} d(v)v^{s-1}\Delta(1/v) = 0$, with $\Delta(r) := T_0(i\cdot r^2, \beta, 4w) - \theta_{T,2\beta}(r)$. $|\varphi(r)| \leq Kg(r, j, 4(w + 1))$, with $j = \max\{2\mathrm{Re}(\alpha), 4w\}$, for each of $\varphi(r) = T_0(i\cdot r^2, \beta, 4w)$, $\varphi(r) = \theta_{T,2\beta}(r)$, and so also for $\varphi(r) = \Delta(r)$. Select an x with $j < x < 4(w + 1)$. $\int_{v>0} d(v)v^{x-1}\Delta(1/v) = 0$ implies $\Delta(r) = 0$ for almost all r > 0. Δ(r) is continuous. Hence $T_0(i\cdot r, \beta, 4w) = \sum_{n\geq 1} (n^{-2\beta})P_{4w}(\pi\cdot r/(n^2), \beta)$.

Say $\beta \geq 0$ and $P_{4w}(v, \beta) > 0$ for $v > 0$. Then $T_0(i\cdot r, \beta, 4w) > 0$.

**Proof of (2).** Say $\beta_1 \geq 0$ and $w \geq 0$. As in the proof of (2) of the previous corollary: $T_0(i\cdot r, \beta, 4w) = h_0^{<\alpha>}(r)$, with $h_0(r) = \pi^{-1+\beta}\cdot½H(\pi r, \beta, 2w)$. Assume $\beta \geq 0$. $H(v, \beta, 2w) > 0$, for all $v > 0$, by Lemma 5.4 (2). Then Lemma 3.2 (4) gives $T_0(i\cdot r, \beta, 4w) > 0$, for $r > 0$.

**Lemma 6.2 Positivity** *Suppose $\beta \geq 0$.*
(1) *Assume $w \geq 0$. $P_{4w}(v, \beta)$ strictly increases with v on $[0, \pi]$.*
(2) *Assume $w \geq 1$. $P_{4w}(v, \beta)$ is positive, for all $v > 0$.*

**Proof of (1).** Say $0 < r \leq 1$. $P_{4w}(\pi r, \beta) = \pi^{-1+\beta} \sum_k \pi^{2k}(a_k'\cdot u^k - \pi^2\cdot a_{k+1}'\cdot u^{k+1})$. Here $k \geq w + 1$ and $k \equiv (w + 1) \bmod(2)$. $u = r^2$. Also $a_k' = a_k/b_k$, with $a_k = 1/\Gamma(1 + \beta + 2k)$ and $b_k = (2k - \alpha)\zeta(2\beta + 4k)$, where $\alpha = ½ - \beta$.

$(d/du)(a_k'\cdot u^k - \pi^2\cdot a_{k+1}'\cdot u^{k+1}) = (u^{k-1}(k + 1)\cdot a_{k+1}'/\zeta(2\beta + 4k))t(u, k, \beta)$, with $t(u, k, \beta) = t(k, \beta)\zeta(2\beta + 4(k + 1)) - \pi^2\zeta(2\beta + 4k)u$, where $t(k, \beta) = (1 + \beta + 2k)2k(1 + \beta/(2(k + 1)))(1 + 2/(2k - ½ + \beta))$. Say $\beta \geq 0$ and $k \geq 1$. $t(u, k, \beta) \geq t(1, k, \beta) > t(k, \beta) - \pi^2\zeta(4)$.

Say $k \geq 2$. $t(k, \beta) \geq (1 + 2\cdot2)2\cdot2$. Now $\pi^2 < 10$ (since $\pi < 3 + 1/7$). $\zeta(4) < 2$. So $t(k, \beta) - \pi^2\zeta(4) > 0$. Then $t(u, k, \beta) > 0$.

Consider $t(1, \beta) - \pi^2\zeta(4)$. $t(1, \beta) = (7 + 2\beta)(1 + 3/(3 + 2\beta))(1 + \beta/4)$. $\zeta(4) = \pi^4/90$. $\pi^2\zeta(4) < 10^2/9$.

Claim *Say $\beta \geq 0$. $t(1, \beta) \geq 11 + 1/5$.*

Proof of Claim If $0 \leq \beta \leq 1$, then $t(1, \beta) \geq 7(1 + 3/5)$. If $\beta \geq 1$, then $t(1, \beta) \geq (7 + 2)(1 + 1/4)$.

Thus when $\beta \geq 0$, $t(1, \beta) - \pi^2\zeta(4) > 0$. Then $t(u, 1, \beta) > 0$.

**Proof of (2).** Say $w \geq 1$. The previous Unconditional theorem 6.1 gives $(-1)^w f(s, \beta) = \int_{v>0} d(v)v^{s-1}P_{4w}(\pi\cdot v^{-2}, \beta)$, when $\beta_1 := \mathrm{Re}(\beta) \geq 0$ and $4w < \mathrm{Re}(s) < 4(w + 1)$, with $P_{4w}(\pi\cdot z, \beta) = \sum_{n\geq1} (\mu(n)/(n^{2\beta}))T_0(iz/(n^2), \beta, 4w)$. According to Lemma 6.1 (2), $T_0(ir, \beta, 4w) > 0$, for all $r > 0$, provided $\beta \geq 0$. However $\mu(n)$ is -1 when n is the product of an odd number of distinct primes. The method employed next "averages out" those -1's.

Fix $\beta$.

$$f(s, \beta) = (1/\zeta(2\beta + s))(1/(s - 2\alpha))\cdot\pi^{-1+\beta}(½)\cdot\pi^{s/2}\cdot F(s/2, \beta),$$

with $\alpha = ½ - \beta$. Say $w \geq 1$. Set $m = 2w$. The translation relation for $F(z, \beta)$ yields

$$(-1)^w f(S + 2m, \beta) = \pi^m(1/(S + 2m - 2\alpha))j(2\beta + S, m)f_0(S, \beta).$$

(See §3.)

Assume $Re(S) > -2(1 + \beta_1)$. Apply Lemma 3.3. The multiplication of $1/\zeta(2\beta + S + 2m)$ by $1/((1 + \beta + S/2)_m)$ yields $j(2\beta + S, m)$, which has a Mellin transform representation with a positive density when $\beta$ is real. $j(2\beta + S, m) = \int_{0 < v \leq 1} (dv)v^{S-1}g_2(v, m)$, with $g_2(v, m) := v^{2\beta}E(v, m)$. $0 < E(v, m) < K(m)v^2$, when $0 < v < 1$.

Assume $Re(S) > 2(-m + ½ - \beta_1)$. Then $1/(S + 2m - 2\alpha) = \int_{0 < v \leq 1} (dv)v^{S-1}g_1(v, m)$, with $g_1(v, m) := v^{2m - 2\alpha}$. $g_1(v, m) > 0$, for all $v > 0$, when $\beta$ is real.

Assume $\beta_1 \geq 0$ and $0 < Re(S) < 4$.

Apply the Corollary 6.1 (1).

$$f_0(S, \beta) = \int_{v > 0} d(v)v^{S-1}\cdot\pi^{-1+\beta}H(\pi\cdot v^{-2}, \beta),$$

with $q_0(H) \geq 2$ and $j_0(H) \leq 0$. The last integral converges absolutely.

Then

$$(-1)^w f(S + 2m, \beta) = \int_{v > 0} d(v)v^{S-1}G(v, m),$$

where

$$G(v, m) := \pi^m\iint_D d(v_1)d(v_2)v_1^{-1}g_1(v, m)v_2^{-1}g_2(v, m)\cdot\pi^{-1+\beta}H(\pi\cdot(v_1v_2/v)^2, \beta),$$

with $D := (0,1]^2$. The last two integrals converge absolutely. $G(v, m)$ is continuous in the $v > 0$.

Say $4w < Re(s) < 4(w + 1)$. Set $m = 2w$. $(-1)^w f(s, \beta) = \int_{v > 0} d(v)v^{s-1}v^{-2m}G(v, m)$. The previous Unconditional theorem 6.1 (2) gives

$$(-1)^w f(s, \beta) = \int_{v > 0} d(v)v^{s-1}P_{4w}(\pi\cdot v^{-2}, \beta).$$

The last integral converges absolutely. $P_{4w}(r, \beta)$ is continuous in the $r > 0$. Hence $P_{4w}(\pi\cdot v^{-2}, \beta) = v^{-2m}G(v, m)$.

Assume $\beta \geq 0$. G(v, m) is positive for all v > 0, since that holds for each of the functions $g_1(v, m)$, $g_2(v, m)$ and $H(v, \beta)$. So $P_{4w}(r, \beta) > 0$, for r > 0.

Fix $\beta \geq 0$. f(s, β) is a meromorphic characteristic function on the complex domain of s with $Re(s) \geq 4$.

**Results when β = ¼ .**

Take β = ¼ in Lemma 6.2 to obtain the following Theorem 3.3 and Lemma 3.2 presented in Part I, §3.

**Theorem 3.3** *$P_0(v)$ is strictly increasing with v for $0 \leq v \leq \pi$ . $P_0(0) = 0$.*

**Lemma 3.2** *Fix $w \geq 1$.*
(i) *$P_{4w}(v)$ is strictly monotone increasing in v on the interval $0 \leq v \leq \pi$. $P_{4w}(0) = 0$.*
(ii) *$P_{4w}(v) > 0$ for v > 0.*

Main unconditional theorem (4) (i) and Lemma 3.2 (ii) together give the following theorem stated in Part I, §3.

**Main unconditional theorem**
**(4)** ***Say $w \geq 1$.***
**(i)** ***On $V_{4w}$, f(s) is analytic and***

$$(-1)^w f(s) = \int_R d(y) e^{sy} P_{4w}(\pi e^{-2y}),$$

***with $P_{4w}(z)$ entire.***
**(ii)** ***$P_{4w}(\pi e^{-2y}) > 0$ for all real y.***

f(-s) = f(s) implies that for $w \leq -1$: $P_{4w}(\pi v) = P_{-4(w+1)}(\pi/v)$, for v > 0. Thus f(s) is a meromorphic characteristic function on the domain of s with $|Re(s)| \geq 4$.

**Metric norms and analytic characteristic functions.**

Fix $\beta \geq 0$ and $w \geq 1$. The positivity of the transform density $P_{4w}(v, \beta)$, when v > 0, has a geometric consequence. That involves the internal metric norms obtained by the author in Part VI (see A. Csizmazia [7]), as a generalization of the Fourier metric norms of J. von Neumann and I. J. Schoenberg [11]. The metric arising from f(s, β) is presented next.

**Definitions of m(t), $d(t_1, t_2)$.** Fix $\beta \geq 0$. Say $x > 4$ and x is not a multiple of four. Fix x. Let t be real. $m(t) := |1 - (n(x, \beta)/n(x + it, \beta))|^{1/2}$. $d(t_1, t_2) := m(t_1 - t_2)$, for real $t_1, t_2$.

The following is a corollary of Unconditional theorem 6.1 (2) and Lemma 6.2 (2) herein together with Corollary 2.2 of Part VI.

**Corollary 6.2** *Say $\beta \geq 0$, $x > 4$, and x is not a multiple of four.*

(1) *The associated m(t) is a metric norm in t on the real line.*

(2) *The associated $d(t_1, t_2)$ is a (finite-valued) translation invariant metric in real $t_1, t_2$.*

**Metric result when β = ¼ .**

Assume $x > 4$ or $x < -4$, with x is not a multiple of four. Fix x. Say t is real. Set $m(t) := |1 - (n(x)/n(x + it))|^{1/2}$. Let $d(t_1, t_2) := m(t_1 - t_2)$ for real $t_1, t_2$. In that case (1) and (2) of the previous corollary hold, provided Corollary 2.2 of Part VI does. See Part I, §3, A geometric consequence of the Main unconditional theorem (4).

Corollary 6.2 with β = ¼ has a conditional generalization to all real x other than the multiples of four. That result is established in Part VI. See Metric norms and analytic characteristic functions of Part I, §7 and Part IV, §4.

**References**

[1] L. Ahlfors, *Complex Analysis,* 3rd Ed., McGraw–Hill, New York, (1978).

[2] G. Andrews, R. Askey, R. Roy, *Special Functions,* Encyclopedia of Mathematics, Cambridge University Press, pp 18-22, (1999).

[3] T. M. Apostol, *Introduction to Analytic Number Theory*, Springer Verlag, New York, (1976).

[4] H. Buchholz, *The Confluent Hypergeometric Function with Special Emphasis on its Applications*, Springer-Verlag, New York, (1969).

[5] A. Csizmazia, On the Riemann zeta-function, Part I: Outline

[6] A. Csizmazia, On the Riemann zeta-function and meromorphic characteristic functions. On the Riemann zeta-function, Part IV.

[7] A. Csizmazia, A generalization of the Fourier metric geometries of J. von Neumann and I. J. Schoenberg. On the Riemann zeta-function, Part VI.

[8] A. Csizmazia, On the Riemann zeta-function: A relation of its nonreal zeros and first derivatives thereat to its values on ½ + 4**N**. Part V.

[9] Erdélyi et al., *Higher transcendental functions,* vol; 1, ch. 6, On Kummer Functions, McGraw-Hill Book Co., Inc., New York, (1953).

[10] E. Lukacs, *Characteristic Functions,* 2nd ed., Griffin, London, (1970).

[11] J. von Neumann and I. J. Schoenberg (1941), Fourier integrals and metric geometry, *Trans. Amer. Math. Soc.* **50**, 226–251.

[12] L. J. Slater, Confluent Hypergeometric Functions, ch. 13, *Handbook of Mathematical Functions*, M. Abramowitz and I. A. Stegun, Editors, National Bureau of Standards, Applied Mathematics Series - 55 Tenth Printing, (1972). http://www.math.sfu.ca/~cbm/aands/intro.htm#006

[13] L. J. Slater, *Confluent Hypergeometric Functions,* Cambridge, England: Cambridge University Press, (1960).

[14] F. G. Tricomi, *Fonctions hypergéométriques confluentes.* Gauthier-Villars, Paris, (1960).

[15] Eric W. Weisstein, Generalized Hypergeometric Function, from *Math World* - A Wolfram Web Resource. http://mathworld.wolfram.com/GeneralizedHypergeometricFunction.html